\documentclass[11pt]{article}
\usepackage{latexsym}
\usepackage{amscd}
\usepackage{amssymb}
\usepackage{amsmath}
\usepackage{amsthm}
\usepackage{enumerate}
\usepackage{verbatim}
\def\co{\colon\thinspace}
\newtheorem{thm}{Theorem}[section]
\newtheorem{cor}[thm]{Corollary}

\newtheorem{lem}[thm]{Lemma}
\newtheorem{prop}[thm]{Proposition}
\newtheorem{dfn}[thm]{Definition}

\newtheorem{Example}[thm]{Example}
\newenvironment{ex}{\begin{Example}\rm}{\end{Example}}
\newtheorem{remark}[thm]{Remark}
\newenvironment{rmk}{\begin{remark}\rm}{\end{remark}}
\newtheorem{Fact}[thm]{Fact}

\begin{document}
\abovedisplayskip=6pt plus3pt minus3pt
\belowdisplayskip=6pt plus3pt minus3pt
\title{Pinching, Pontrjagin classes, and \\
negatively curved vector bundles}
\author{Igor Belegradek}
\date{}
\maketitle
\begin{abstract} 
We prove several finiteness results for \@
the class $\mathcal M_{a,b,\pi, n}$ 
of $n$-manifolds that have fundamental groups
isomorphic to $\pi$ and that can be given complete Riemannian 
metrics of sectional curvatures within $[a,b]$ where $a\le b<0$.
In particular,
if $M$ is a closed negatively curved manifold of dimension at least three,
then only finitely many manifolds in the class $\mathcal M_{a,b,\pi_1(M), n}$ 
are total spaces of vector bundles over $M$.
Furthermore,
given a word-hyperbolic group $\pi$ and an integer
$n$ there exists a positive $\epsilon=\epsilon(n,\pi)$ 
such that the tangent bundle of any manifold in the class
$\mathcal M_{-1-\epsilon, -1, \pi, n}$ has
zero rational Pontrjagin classes.
\end{abstract}
\section{Introduction}
\label{Introduction}
According to the Cartan-Hadamard theorem, the universal cover
of any complete negatively curved manifold is diffeomorphic
to the Euclidean space.
Surprisingly, beyond this fact little is known
about topology of {\it infinite}
volume complete negatively curved manifolds. 
For example, so far there has been found 
no restriction on the fundamental groups of such manifolds
except being the fundamental groups of
aspherical manifolds.

In this paper we study topology of pinched negatively curved manifolds.
(We call a manifold {\it pinched negatively
curved} if it admits a complete Riemannian metric
with sectional curvatures
bounded between two negative constants.)
For instance, the fundamental groups of 
pinched negatively curved manifolds
have the property that any amenable subgroup must be 
finitely generated and virtually nilpotent~\cite{BS, Bow2}.
Examples include closed negatively curved manifolds
and complete locally symmetric negatively curved manifolds.
Also, elementary warped product construction gives
pinched negatively curved metrics on 
the direct products of pinched negatively
curved manifolds with Euclidean spaces~\cite{FJ5}.
Furthermore, the total space of any vector bundle over
a closed negatively curved manifold is pinched negatively
curved~\cite{And}.

According to a general pinching principle, 
a negatively curved manifold with sectional curvatures
close to $-1$ ought to be topologically similar to a hyperbolic manifold. 
In its strongest form this principle fails even for closed manifolds.
In fact, for each $n\ge 4$ and any $\epsilon>0$, 
M.~Gromov and W.~Thurston~\cite{GT} 
found a closed negatively curved $n$-manifold with sectional
curvatures within $[-1-\epsilon, -1]$ that is not diffeomorphic to
a hyperbolic manifold.
(Any closed negatively
curved $3$-manifold is diffeomorphic to
a hyperbolic one if Thurston's hyperbolization conjecture is true.)
However, the pinching principle holds
if the pinching constant $\epsilon$ is allowed to depend
on the fundamental group, namely,
there exists an $\epsilon=\epsilon(\pi) >0$ such that
any closed manifold in the class $\mathcal M_{-1-\epsilon, -1, \pi, n}$
is diffeomorphic to a hyperbolic manifold~\cite{Bel4}.
In even dimensions Gromov~\cite{Gro4} proved a stronger 
pinching theorem: any closed even-dimensional Riemannian 
manifold with sectional curvatures
within $[-1-\epsilon, -1]$ is diffeomorphic to a hyperbolic manifold
where $\epsilon$ depends on the dimension and the Euler characteristic
of the manifold.

The following theorem restricts the topology
of strongly pinched negatively curved manifolds without assuming
compactness.

\begin{thm}\label{intro:pinching}
Let $\pi$ be the fundamental group of a finite aspherical cell complex. 
Suppose that $\pi$ is not virtually nilpotent
and that $\pi$ does not split as a 
nontrivial amalgamated product
or an HNN-extension over a virtually nilpotent group.
Then, for any positive integer $n$, there exists
an $\epsilon=\epsilon(\pi, n)>0$ such that any manifold in
the class $\mathcal M_{-1-\epsilon, -1, \pi, n}$ 
is tangentially homotopy equivalent to
an $n$-manifold of constant negative curvature.
\end{thm}

Recall that a homotopy equivalence of manifolds
$f:M\to N$ is called tangential if the vector bundles
$f^\#TN$ and $TM$ are stably isomorphic.
If we weaken the assumption ``$\pi$ is the fundamental group of a 
finite aspherical cell complex'' to 
``$\pi$ is finitely presented'', then 
the same conclusion holds without the word ``tangentially''.

Any manifold of constant sectional curvature
has zero rational Pontrjagin classes \cite{Ave},
hence, applying the Mayer-Vietoris sequence and
the accessibility result of Delzant
and Potyagailo~\cite{DP} (reviewed in~\ref{s:Accessibility}),
we deduce the following.

\begin{thm}
Let $\pi$ be a finitely presented group with finite
$4k$th Betti numbers for all $k>0$.
Assume that any nilpotent subgroup of $\pi$ has 
cohomological dimension $\le 2$.
Then for any $n$ there exist $\epsilon=\epsilon(n,\pi)>0$ 
such that the tangent bundle of any manifold in the class
$\mathcal M_{-1-\epsilon, -1, \pi, n}$ 
has zero rational Pontrjagin classes.
\end{thm}

In particular, the following is true
because any nilpotent subgroup of
a word-hyperbolic group has cohomological dimension $\le 1$.

\begin{cor}
Let $\pi$ be a word-hyperbolic group.
Then for any $n$ there exist $\epsilon=\epsilon(n,\pi)>0$ 
such that the tangent bundle of any manifold in the class
$\mathcal M_{-1-\epsilon, -1, \pi, n}$ 
has zero rational Pontrjagin classes.
\end{cor}

As we showed in~\cite{Bel5}, the class 
$\mathcal M_{a, b, \pi, n}$ falls into 
finitely many tangentially homotopy types under 
mild assumptions on the group $\pi$.
However,~\cite{Bel5} does not provide explicit bounds
on the number of tangentially homotopy inequivalent
manifolds in the class $\mathcal M_{a, b, \pi, n}$.
The only exception is the class 
$\mathcal M_{-1-\epsilon, -1, \pi, n}$
in the theorem~\ref{intro:pinching}
because the number of tangentially homotopy inequivalent
real hyperbolic manifolds is bounded by
the number of connected components of the representation variety 
$\text{Hom}(\pi,\text{Isom}(\mathbf{H}^n_{\mathbb R}))$~\cite{Bel3}.
The latter can be estimated in terms of $n$ and the numbers
of generators and relators of $\pi$.

M.~Anderson~\cite{And} showed that the total space 
of any vector bundle over a closed negatively curved manifold $M$ 
can be given a complete Riemannian metric of sectional curvature 
within $[a,b]$ for some $a\le b<0$. 
However, the theorem below says that, once pinching $a/b$ is fixed, 
total spaces of only finitely many vector bundles over $M$ of a given
rank admit metrics of sectional curvature within $[a,b]$.
Note that the set of isomorphism classes
of vector bundles of the same rank over $M$ is {\it infinite} 
provided certain Betti numbers of $M$ are nonzero 
(see~\ref{vector bundles infinite}).

\begin{thm}\label{intro:vector bundles}
Let $M$ be a closed negatively curved
manifold with $\dim(M)\ge 3$ and let $k>0$
be a positive integer and $a\le b<0$ be real numbers. 
Then, up to isomorphism, there exist only finitely many
rank $k$ vector bundles over $M$ whose total spaces 
admit complete Riemannian metrics with sectional
curvatures within $[a,b]$.
\end{thm}

As we explain in the appendix, 
the isomorphism type of a vector bundle $\xi$ over a finite 
cell complex is determined, up to finitely many possibilities,
by its Euler class and Pontrjagin classes
(for non-orientable bundles one looks at the Euler class of
the orientable two-fold-pullback).
Pontrjagin classes depend on the tangent bundle of 
the total space $E(\xi)$ while the Euler class
can be computed via intersections in $E(\xi)$.
Finiteness results for tangent bundles of pinched
negatively curved manifolds were obtained in~\cite{Bel5};
in the present paper we prove similar results for
intersections.

A homotopy equivalence $f:N\to L$ of orientable $n$-manifolds is called
{\it intersection preserving} if, 
for some orientations of $N$ and $L$, the intersection number of
any pair of homology classes of $N$ is equal to
the intersection number of their $f$-images in $L$.  
More generally, a homotopy equivalence $f$ of
nonorientable manifolds is called {\it intersection preserving}
if $f$ preserves first Stiefel-Whitney classes and 
the lift of $f$ to
orientable two-fold covers is an intersection preserving
homotopy equivalence. 
For example, if $f$ is homotopic to a 
homeomorphism, then $f$ is an intersection preserving
homotopy equivalence. 

\begin{thm}\label{intro:no split means bound on intersections}
Let $\pi$ be a finitely presented group with finite 
Betti numbers.
Assume that $\pi$ does not split 
as a nontrivial amalgamated product or
an HNN-extension over a virtually nilpotent group.
Then, for any $a\le b<0$,
the class $\mathcal M_{a, b, \pi, n}$ breaks into finitely
many intersection preserving
homotopy types.
\end{thm}

Theorem~\ref{intro:no split means bound on intersections} was first
proved in~\cite{Bel3} for oriented locally symmetric negatively 
curved manifolds. 
This special case 
is somewhat easier to handle due to the fact that
theory of convergence discussed in
the section~\ref{Convergence, splittings and accessibility}
has been thoroughly studied in the constant negative curvature case.
Note that if $a=b=-1$, one can naturally identify 
$\mathcal M_{a,b,\pi, n}$ with the set of injective 
discrete representations of $\pi$ into the isometry group of the
real hyperbolic $n$-space which is a standard object
in Kleinian group theory.
\paragraph{Synopsis of the paper.}
The second section contains background in convergence
of Riemannian manifolds as well as
some results on splittings and accessibility over virtually
nilpotent groups.
The main technical result is proved in the section three
after discussing some invariants of maps and proper
discontinuous actions. 
The forth section is devoted to applications.
Section five is a discussion of certain natural 
pinching invariants.
A bundle-theoretic result is proved in the appendix.
\paragraph{Acknowledgments.}
I am grateful to Werner Ballmann, Mladen Bestvina, 
Thomas Delzant, Martin J.~Dunwoody, Karsten Grove, Lowell E.~Jones,
Misha Kapovich, Vitali Kapovitch, Bruce Kleiner,
Bernhard Leeb, John J.~Millson,
Igor Mineyev, Sergei P.~Novikov, 
Fr\'ed\'eric Paulin, Conrad Plaut,
Jonathan M.~Rosenberg, James A.~Schafer, Harish Seshadri, and
Shmuel Weinberger for helpful discussions and communications.
Special thanks are due to my advisor 
Bill Goldman for his constant interest and support.
This work is a part of my thesis at the University
of Maryland, College Park.
Also I wish to thank Heinz Helling and 
SFB-343 at the University of Bielefeld
for support and hospitality.

\section{Convergence, splittings and accessibility}
\label{Convergence, splittings and accessibility}
The exposition in this section is a variation of
the one given in~\cite{Bel5}. 
Several straightforward lemmas are only stated and 
refered to~\cite{Bel5} for proofs.

By an {\it action} of an abstract group $\pi$ on a space $X$
we mean a group homomorphism $\rho:\pi\to\mathrm{Homeo}(X)$.
An action $\rho$ is called {\it free} if $\rho(\gamma)(x)\neq x$
for all $x\in X$ and all $\gamma\in \pi\setminus\text{id}$.
In particular, if $\rho$ is a free action, then $\rho$ is injective. 

\subsection{Equivariant pointed Lipschitz topology}
\label{Equivariant pointed Lipschitz topology}
Let $\Gamma_k$ be a discrete
subgroup of the isometry group
of a complete Riemannian manifold $X_k$
and $p_k$ be a point of $X_k$.
The class of all such triples $\{(X_k, p_k, \Gamma_k)\}$ 
can be given the so-called 
equivariant pointed Lipschitz topology~\cite{Fuk}; 
when $\Gamma_k$ is
trivial this reduces to the usual pointed Lipschitz topology.
For convenience of the reader
we give here some definitions borrowed from~\cite{Fuk}.

For a group $\Gamma$ acting on a pointed metric space $(X,p,d)$
the set $\{\gamma\in\Gamma: d(p, \gamma(p))<r\}$ is denoted by $\Gamma(r)$.
An open ball in $X$ of radius $r$ with center at $p$ is denoted by
$B_r(p,X)$.

For $i=1,2$, let $(X_i, p_i)$ be a pointed
complete metric space with the distance function
$d_i$ and let $\Gamma_i$ be a discrete group of isometries
of $X_i$. In addition, assume that $X_i$ is a $C^\infty$--manifold.
Take any $\epsilon>0$.

Then a quadruple $(f_1, f_2, \phi_1, \phi_2)$ of maps 
$f_i:B_{1/\epsilon}(p_i, X_i)\to B_{1/\epsilon}(p_{3-i}, X_{3-i})$
and $\phi_i:\Gamma_i(1/3\epsilon)\to\Gamma_{3-i}$
is called an $\epsilon$--{\it Lipschitz approximation}
between the triples
$(X_1,p_1,\Gamma_1)$ and $(X_2,p_2,\Gamma_2)$
if the following seven condition hold:

$\bullet$ $f_i$ is a diffeomorphism onto its image;

$\bullet$ for each $x_i\in B_{1/3\epsilon}(p_i, X_i)$
and every $\gamma_i\in\Gamma_i(1/3\epsilon)$,  
$f_i(\gamma_i(x_i))=\phi_i(\gamma_i)(f_i(x_i))$;

$\bullet$ for every $x_i, x_i^\prime\in B_{1/\epsilon}(p_i, X_i)$,
$e^{-\epsilon}<
d_{3-i}(f_i(x_i),f_i(x_i^\prime))/d_i(x_i,x_i^\prime)<e^\epsilon$; 

$\bullet$ $f_i(B_{1/\epsilon}(p_i, X_i))\supset
B_{(1/\epsilon)-\epsilon}(p_{3-i}, X_{3-i})$ and 
$\phi_i(\Gamma_i(1/3\epsilon))\supset\Gamma_{3-i}(1/3\epsilon-\epsilon)$;

$\bullet$ $f_i(B_{(1/\epsilon)-\epsilon}(p_i, X_i))\supset
B_{1/\epsilon}(p_{3-i}, X_{3-i})$ and 
$\phi_i(\Gamma_i(1/3\epsilon-\epsilon))\supset\Gamma_{3-i}(1/3\epsilon)$;

$\bullet$ $f_{3-i}\circ f_i|_{B_{(1/\epsilon)-\epsilon}(p_i, X_i)}=\mathrm{id}$
and $\phi_{3-i}\circ\phi_i|_{\Gamma_i(1/3\epsilon-\epsilon)}=\mathrm{id}$;

$\bullet$ $d_{3-i}(f_i(p_i),p_{3-i})<\epsilon$.

We say a sequence of triples $(X_k, p_k, \Gamma_k)$
{\it converges} to $(X, p, \Gamma)$ in the 
equivariant pointed Lipschitz topology
if for any $\epsilon>0$
there is $k(\epsilon)$ such that for all $k>k(\epsilon)$,
there exists an $\epsilon$--Lipschitz approximation
between $(X_k, p_k, \Gamma_k)$ and $(X, p, \Gamma)$.
If all the groups $\Gamma_k$ are trivial,
then $\Gamma$ is trivial; in this case 
we say that that $(X_k,p_k)$ converges to $(X,p)$ in the pointed
Lipschitz topology. 

Note that if $X_k$ is a complete
Riemannian manifold for all $k$,
then the space $X$ is necessarily a $C^\infty$--manifold with 
a complete $C^{1,\alpha}$--Riemannian metric~\cite{GW}. 
If each $X_k$ is a Hadamard manifold,
then for any $x_k\in X_k$, the sequence $(X_k,x_k)$
is precompact in the pointed Lipschitz topology
because the injectivity radius of $X_k$ at $x_k$
is uniformly bounded away from zero~\cite[p.132]{Fuk}.

\begin{rmk}
Equivariant pointed Lipschitz topology 
is closely related to the so-called
Chabauty topology used 
in Kleinian group theory~\cite{CEG},~\cite{BP}.
Indeed, let $\Gamma_k$ be a sequence of discrete
subgroups of the isometry group of a complete Riemannian manifold
$X$ (e.g.~a hyperbolic space).
Then $\Gamma_k$ converges in the 
Chabauty topology to a {\it discrete}
group $\Gamma$ if and only if for each $p\in X$
$(X,\Gamma_k, p)$ converges to $(X,\Gamma, p)$
in the equivariant pointed Lipschitz topology.
\end{rmk}   

\subsection{Pointwise convergence topology}
\label{Pointwise convergence topology}
Suppose that, for some $p_k\in X_k$, 
the sequence  $(X_k, p_k)$ converges to $(X,p)$
in the pointed Lipschitz topology, 
i.e.,\! for any $\epsilon>0$
there is $k(\epsilon)$ such that for all $k>k(\epsilon)$,
there exists an $\epsilon$--Lipschitz approximation $(f_k, g_k)$
between $(X_k, p_k)$ and $(X, p)$.
We say that a sequence $x_k\in X_k$ {\it converges} to
$x\in X$ if for some $\epsilon$
$$d(f_k(x_k), x)\to 0\ \ \text{as}\ \ k\to\infty$$
where $d(\cdot,\cdot)$ is the distance function on $X$
and $f_k$ comes from the $\epsilon$--Lipschitz approximation
$(f_k, g_k)$ between $(X_k, p_k)$ and $(X, p)$.
Trivial examples: if $(X_k, p_k)$ converges to $(X,p)$
in the pointed Lipschitz topology, then $p_k$ converges to $p$;
furthermore, if $x\in X$, the sequence $g_k(x)$ converges to $x$.

Given a sequence of isometries $\gamma_k\in\mathrm{Isom}(X_k)$
we say that $\gamma_k$ {\it converges}, if for any $x\in X$ and any 
sequence $x_k\in X_k$ that converges to $x$,
$\gamma_k(x_k)$ converges. 
The limiting transformation $\gamma$ that takes $x$ to
the limit of $\gamma_k(x_k)$
is necessarily an isometry of $X$.
Furthermore, if $\gamma_k$ and $\gamma_k^\prime$
converge to $\gamma$ and $\gamma^\prime$ respectively,
then $\gamma_k\cdot\gamma_k^\prime$ converges
to $\gamma\cdot\gamma^\prime$.
In particular, $\gamma_k^{-1}$ converges to $\gamma^{-1}$
since the identity maps $\text{id}_k:X_k\to X_k$
converge to $\text{id}:X\to X$.

Let $\rho_k:\pi\to\mathrm{Isom}(X_k)$ 
be a sequence of isometric actions of a group $\pi$ on $X_k$.
We say that a sequence of actions $(X_k, p_k, \rho_k)$ 
{\it converges in the pointwise 
convergence topology} if $\rho_k(\gamma)$ converges 
for every $\gamma\in\pi$.
The limiting map $\rho:\Gamma\to\text{Isom}(X)$ that takes 
$\gamma$ to the limit of $\rho_k(\gamma)$ is necessarily
a homomorphism. 
If $\pi$ is generated by a finite set $S$, then
in order to prove that $\rho_k$ converges in the pointwise
convergence topology it suffices to check that 
$\rho_k(\gamma)$ converges, for every $\gamma\in S$.   
 
It is worth clarifying that the term ``pointwise convergence'' refers
to the convergence of group action rather than
individual isometries. 
The definitions are set up so that
individual isometries converge 
``uniformly on compact subsets''.
The motivation comes from the following example.

\begin{ex}
Let $X$ be a complete Riemannian manifold (e.g.~a hyperbolic space).
Consider the isometry group
$\mathrm{Isom}(X)$ equipped with
compact-open topology and let $\pi$ be a group. 
The space $\mathrm{Hom}(\pi,\mathrm{Isom}(X))$ has a natural
topology (which is usually called ``algebraic topology''
or ``pointwise convergence topology''), namely
$\rho_k$ is said to converge to $\rho$ if, for each $\gamma\in\pi$,
$\rho_k(\gamma)$ converges to $\rho(\gamma)$ in the Lie group 
$\mathrm{Isom}(X)$.  
Note that if $\pi$ is finitely generated, this topology
on $\mathrm{Hom}(\pi,\mathrm{Isom}(X))$
coincide with the compact-open topology.
Certainly, for any $p\in X$, 
the constant sequence $(X,p)$ converges to itself
in pointed Lipschitz topology. 
Then the sequence $(X, p, \rho_k)$
converges in the pointwise convergence topology
(as defined in this section) 
if and only if 
$\rho_k\in\mathrm{Hom}(\pi,\mathrm{Isom}(X))$
converges in the algebraic topology.
(Indeed, $\rho_k(\gamma)$ converges to $\rho(\gamma)$ in 
$\mathrm{Isom}(X)$ iff $\rho_k(\gamma)$ converges to $\rho(\gamma)$
uniformly on compact subsets. 
In particular, the latter implies that
$\rho_k(\gamma)(x_k)\to\rho(\gamma)(x)$ for any $x_k\to x$. 
Conversely, if $\rho_k(\gamma)(x_k)\to\rho(\gamma)(x)$ for any $x_k\to x$,
then $\rho_k(\gamma)$ converges to $\rho(\gamma)$
uniformly on compact subsets~\cite[4.7, Lemma~5]{KN}.)
\end{ex}

A sequence of actions $(X_k, p_k, \rho_k)$
is called {\it precompact in the pointwise convergence
topology} if every subsequence of
$(X_k, p_k, \rho_k)$ has a subsequence that
converges in the pointwise
convergence topology.

Repeating the proof of~\cite[4.7]{KN}, 
it is easy to check that 
a sequence of isometries $\gamma_k\in\text{Isom}(X_k)$ 
has a converging subsequence if, 
for some converging sequence $x_k\in X_k$, 
the sequence $d_k(x_k,\gamma_k(x_k))$
is bounded 
(where $d_k(\cdot,\cdot)$ is the distance function on $X_k$).

Suppose that $\pi$ is a countable (e.g. finitely generated)
group and assume that
for each $\gamma\in\pi$
the sequence $d_k(p_k,\rho_k(\gamma)(p_k))$
is bounded. Then $(X_k, p_k, \rho_k)$
is precompact in the pointwise convergence topology. 
(Indeed, let $\gamma_1\dots\gamma_n\dots$ be the list of 
all elements of $\pi$. 
Take any subsequence $\rho_{k,0}$ of $\rho_k$. 
Pass to subsequence $\rho_{k,1}$ of $\rho_{k,0}$ so that
$\rho_{k,1}(\gamma_1)$ converges.
Then pass to subsequence $\rho_{k,2}$ of $\rho_{k,1}$
such that $\rho_{k,2}(\gamma_2)$ converges, etc.
Then $\rho_{k,k}(\gamma_n)$ converges for every $n$.)

Note that if $\pi$ is generated by a finite set $S$, then
to prove that $\rho_k$ is precompact
it suffices to check that 
$d_k(p_k,\rho_k(\gamma)(p_k))$ is bounded, 
for all $\gamma\in S$ because it implies that 
$d_k(p_k,\rho_k(\gamma)(p_k))$ is bounded, 
for each $\gamma\in\pi$.

The following two lemmas can
be easily deduced from definitions;
the reader is refered to~\cite{Bel5} for details.

\begin{lem}\label{geom contains alg, etc} 
Let $\rho_k:\pi\to\mathrm{Isom}(X_k)$
be a sequence of isometric actions
of a discrete group $\pi$ on complete Riemannian
$n$-manifolds $X_k$ such that 
$\rho_k(\pi)$ acts freely.
If the sequence
$(X_k, p_k, \rho_k(\pi))$ converges in the 
equivariant pointed Lipschitz topology to $(X, \Gamma, p)$
and $(X_k, p_k, \rho_k)$ converges to $(X, p, \rho)$ 
in the pointwise convergence topology, then 

$(1)$ $\Gamma$ acts freely, and

$(2)$ $\rho(\pi)\subset\Gamma$, and

$(3)$ $\ker(\rho)\subset\ker(\rho_k)$, for all large $k$.
\end{lem}

\begin{lem}\label{approximation}
Let $\rho_k:\pi\to\mathrm{Isom}(X_k)$
be a sequence of isometric actions
of a discrete group $\pi$ on complete Riemannian
$n$-manifolds $X_k$ such that 
$\rho_k(\pi)$ acts freely.
Suppose that the sequence
$(X_k, p_k, \rho_k(\pi))$ converges in the 
equivariant pointed Lipschitz topology to $(X, \Gamma, p)$
and $(X_k, p_k, \rho_k)$ converges to $(X, p, \rho)$ 
in the pointwise convergence topology. 

Then, 
for any $\epsilon>0$ and for any finite subset $S\subset\pi$,
there is $k(\epsilon, S)$ with the property that
for each $k>k(\epsilon, S)$ there exists an
$\epsilon$--Lipschitz approximation 
$(f_k, g_k,\phi_k,\tau_k)$
between $(X_k, p_k, \rho_k(\pi))$
and $(X, p,\Gamma)$ such that 
$\phi_k(\rho_k(\gamma))=\rho(\gamma)$ and 
$\rho_k(\gamma)=\tau_k(\rho(\gamma))$
for every $\gamma\in S$.
\end{lem}

\subsection{Applications of the Margulis' lemma}

The following proposition generalizes
a well-known statement in the Kleinian group theory.
Namely, if $\pi$ is a non-virtually-abelian group
and $\rho_k$ is a sequence of injective discrete
representations of $\pi$ into $\mathbf{PSL}(2,\mathbb C)$
that converges algebraically to a representation
$\rho$, then $\rho$ is injective and discrete.
Moreover, the closure of $\{\rho_k(\pi)\}$
in the Chabauty topology consists of discrete groups.
 
\begin{prop}\label{pointwise implies lipsch}
Let $X_k$ be a sequence of Hadamard manifolds
with sectional curvatures in $[a, b]$ for $a\le b<0$
and let $\pi$ be a finitely generated group that 
is not virtually nilpotent.
Let $\rho_k: \pi\to\mathrm{Isom}(X_k)$ be  
an arbitrary sequence of free and isometric actions 
such that $(X_k,p_k,\rho_k)$ converges 
converges to $(X, p,\rho)$
in the pointwise convergence topology.
Then 
\begin{description}
\item{(i)} the sequence $(X_k, p_k, \rho_k(\pi))$ is precompact
in the equivariant pointed Lipschitz topology, and
\item{(ii)} $\rho$ is a free action, in particular
$\rho$ is injective,
\end{description}
\end{prop}

\begin{proof} 
Choose $r$ so that the open ball $B(p,r)\subset X$
contains $\{\rho(\gamma_1)(p),\dots\rho(\gamma_m)(p)\}$ where
$\{\gamma_1,\dots\gamma_m\}$ generate $\pi$.
Passing to subsequence, we assume that $B(p_k,r)$
contains $\{\rho_k(\gamma_1)(p),\dots\rho_k(\gamma_m)(p)\}$.

Show that, for every $k$, 
there exists $q_k\in B(p_k,r)$ such that
for any $\gamma\in\pi\setminus\{\mathrm{id}\}$,
we have $\rho_k(\gamma)(q_k)\notin B(q_k,\mu_n/2)$ where
$\mu_n$ is the Margulis constant.
Suppose not. 
Then for some $k$, the whole ball $B(p_k,r)$
projects into
the thin part $\{\mathrm{InjRad}<\mu_n/2\}$ 
under the projection $\pi_k:X_k\to X_k/\rho_k(\pi)$.
Thus the ball $B(p_k, r)$ lies in a connected component $W$
of the $\pi_k$--preimage of the thin part of $X_k/\rho_k(\pi)$. 
According to~\cite[p111]{BGS}, the stabilizer of
$W$ in $\rho_k(\pi)$ is virtually nilpotent and,
moreover, the stabilizer contains every element
$\gamma\in\rho_k(\pi)$ with $\gamma(W)\cap W\neq\emptyset$.
Therefore, the whole group $\rho_k(\pi)$ stabilizes $W$.
Hence $\rho_k(\pi)$ must be virtually nilpotent.
As $\rho_k$ is injective, $\pi$ is virtually nilpotent.
A contradiction. 
 
Thus, $(X_k, q_k, \rho_k(\pi))$ is 
Lipschitz precompact~\cite{Fuk} and, hence passing to subsequence,
one can assume that $(X_k, q_k,\rho_k(\pi))$
converges to some $(X, q, \Gamma)$. 

It is a general fact that follows easily from definitions 
that whenever $(X_k, q_k, \Gamma_k)$ converges to $(X, q, \Gamma)$
in the equivariant pointed Lipschitz topology
and a sequence of points $p_k\in X_k$ converges
to $p\in X$, then
$(X_k, p_k, \Gamma_k)$ converges to $(X, p, \Gamma)$
in the equivariant pointed Lipschitz topology.
The proof of $(i)$ is complete.

Pass to a subsequence so that $(X_k, p_k, \rho_k)$
converges to $(X, p, \Gamma)$
in the equivariant pointed Lipschitz topology.
By~\ref{geom contains alg, etc}, $\rho$ is injective. 
Furthermore, $\rho(\pi)$ acts freely because it is
a subgroup of $\Gamma$. Thus, $(ii)$ is proved.
\end{proof}

\subsection{Diverging actions and splittings over virtually
nilpotent groups}

We say that a group $\pi$ splits over a virtually nilpotent 
group if $\pi$
has a nontrivial decomposition into an amalgamated product
or an HNN-extension over a virtually nilpotent group.

Let $\pi$ be a finitely generated group and
let $S\subset\pi$ be a finite
subset that contains $\{\mathrm{id}\}$ and generates $\pi$.
Let $\rho_k: \pi\to\mathrm{Isom}(X_k)$ be  
an arbitrary sequence of free and isometric actions of $\pi$
on Hadamard $n$-manifolds $X_k$.
Assume that the sectional curvatures of $X_k$
lie in $[a, b]$ for $a\le b<0$.

For $x\in X_k$, we denote $D_k(x)$ the diameter
of the set  $\rho_k(S)(x)$.
Set $D_k=\inf_{x\in X_k}D_k(x)$.
Suppose the sequence $D_k$ is bounded. 
Then there exist $x_k\in X_k$ such that 
$D_k(x_k)$ is bounded. Therefore,
as we observed in the section~\ref{Pointwise convergence topology}, 
the sequence $(X_k, x_k,\rho_k)$ 
is precompact in the pointwise convergence topology.
The following lemma shows what happens if $D_k$ is unbounded.

\begin{prop}\label{divergence implies split}
Suppose that the sequence $\{D_k\}$ is unbounded.
Assume that $\pi$ is not
virtually nilpotent.
Then $\pi$ acts on a certain $\mathbb R$-tree without global
fixed points and so that the stabilizer of
any non-degenerate arc is virtually nilpotent.
Furthermore, if $\pi$ is finitely presented, 
then $\pi$ splits over a virtually nilpotent group.
\end{prop}
\begin{proof} 
This proposition is well-known to experts.
First, using work of Bestvina~\cite{Bes} and Paulin~\cite{Pau1, Pau2},
we produce an action of $\pi$
on a real tree and then invoke Rips' machine to get a
splitting over a virtually nilpotent group.
For completeness, we briefly review the argument.

The rescaled pointed Hadamard manifold $\frac{1}{D_k}\!\cdot\! X_k$
has sectional curvature $\le b\cdot D_k\to -\infty$ as $k\to\infty$.
Find $p_k\in X_k$ such that $D_k(p_k)\le D_k+1/k$. 
Consider the sequence of triples 
$(\frac{1}{D_k}\!\cdot\! X_k, p_k, \rho_k)$.
Repeating an argument of Paulin~\cite[\S 4]{Pau2},
we can pass to subsequence that converges 
to a triple $(X_\infty, p_\infty, \rho_\infty)$.
(For the definition of the convergence see~\cite{Pau1, Pau2}.
Paulin calls it ``convergence in the Gromov topology''.) 

The limit space $X_\infty$ is a length space 
of curvature $-\infty$, that is a real tree.
Because of the way we rescaled, 
the limit space has a natural isometric 
action $\rho_\infty$ of $\pi$
with no global fixed point~\cite{Pau1, Pau2}.
Then it is a standard fact that there
exists a unique $\pi$--invariant subtree $T$
of $X_\infty$ that has no proper $\pi$--invariant subtree. 
In fact $T$ is the union of all the axes of all hyperbolic
elements in $\pi$.
Since the sectional curvatures are uniformly bounded
away from zero and $-\infty$, the Margulis lemma implies that
the stabilizer of any 
non-degenerate segment is virtually nilpotent (cf.~\cite{Pau1}). 

Note that any increasing sequence of virtually
nilpotent subgroups of $\pi$ is stationary.
Indeed, since a virtually nilpotent group is amenable, 
the union $U$ of an increasing sequence 
$U_1\subset U_2\subset U_3\subset\dots$
of virtually nilpotent subgroups is also an amenable group.
If the fundamental group of a complete
manifold of pinched negative curvature is amenable,
it must be finitely generated~\cite{BS, Bow2}.
In particular, $U$ is finitely generated,
hence $U_n=U$ for some $n$. 
Thus, the $\pi$--action on the tree $T$ is 
stable~\cite[Proposition 3.2(2)]{BF}.

We summarize that the $\pi$--action
on $T$ is stable, has virtually nilpotent arc stabilizers
and no proper $\pi$--invariant subtree.
Since $\pi$ is finitely presented 
and not virtually nilpotent,
the Rips' machine~\cite[Theorem 9.5]{BF} 
produces a splitting of
$\pi$ over a virtually solvable group.
Any amenable subgroup of $\pi$ must be virtually nilpotent~\cite{BS, Bow2},
hence $\pi$ splits over a virtually nilpotent group.
\end{proof}
%
%
\begin{ex}\label{aspherical do not split}
Let $K$ be a closed aspherical $n$--manifold such that
any nilpotent subgroup of $\pi_1(K)$
has cohomological dimension $\le n-2$.
Then $\pi_1(K)$ does not split as a nontrivial
amalgamated product or HNN-extension over a
virtually nilpotent group~\cite{Bel3}.
In particular, this is the case if $K$ is a closed  aspherical manifold
such that $\dim(K)\ge 3$ and
$\pi_1(K)$ is word-hyperbolic.
\end{ex}

\begin{thm}\label{no split plus pinched means hyperbolic}
Let $\pi$ be a finitely presented group that is not virtually nilpotent
and that does not split over virtually nilpotent subgroups. 
Suppose that $\pi$ is not isomorphic to
a discrete subgroup of the isometry group of the hyperbolic $n$-space.
Then the class 
$\mathcal M_{-1-\epsilon, -1, \pi, n}$ is empty, for some $\epsilon >0$.
\end{thm}
\begin{proof} 
Arguing by contradiction, we assume that
for each $k$ there exists a manifold
$X_k/\rho_k(\pi)\in\mathcal M_{-1-1/k, -1, \pi, n}$
where $\rho_k$ is free, isometric action $\pi$
on the Hadamard manifold $X_k$.

According to~\ref{divergence implies split},
$(X_k, p_k,\rho_k)$ is precompact in the pointwise convergence topology,
for some $p_k\in X_k$.
Pass to subsequence so that $(X_k, p_k,\rho_k)$
converges to $(X, p,\rho)$ in the pointwise convergence topology.
By~\ref{geom contains alg, etc}, 
$\rho$ is a free action of a discrete group $\pi$, in particular
$X/\rho(\pi)$ is a manifold with fundamental 
group isomorphic to $\pi$.

By a standard argument $X$ is an inner metric space
of curvature $\ge -1$ and $\le -1$
in the sense of Alexandrov.
Hence~\cite{Ale} implies that
$X$ is isometric to the real hyperbolic space.
Thus, $\pi$ is the fundamental group of a
hyperbolic manifold $X/\rho(\pi)$ which is a contradiction.
\end{proof}

\subsection{Accessibility over virtually nilpotent groups}
\label{s:Accessibility}
Delzant and Potyagailo have recently proved a powerful
accessibility result which we state below (in torsion-free case)
for reader's convenience.
The definition~\ref{elementary class}
and the theorem~\ref{accessability result} 
are taken from~\cite{DP}.
I am grateful to Thomas Delzant for helpful
discussions.

Recall that a graph of groups is a graph
whose vertices and edges are labeled
with {\it vertex groups} $\pi_v$ and 
{\it edge groups} $\pi_e$ and such that
every pair $(v,e)$ where 
the edge $e$ is incident to the vertex $v$
is labeled with a group monomorphism
$\pi_e\to\pi_v$.
We only consider finite connected graphs of groups.
To each graph of groups one can associate its 
fundamental group which is a result of repeated
amalgamated products and HNN--extensions of vertex groups
over the edge groups (see~\cite{Bau} for more details). 

\begin{dfn}\label{elementary class}
A class $\mathcal E$ of 
subgroups of a torsion free group $\pi$
is called {\it elementary} provided the following four
conditions hold.
\begin{description}
\item[(i)] 
$\mathcal E$ is closed under conjugation in $\pi$;
\item[(ii)] any infinite group from $\mathcal E$ is contained
in a {\bf unique} maximal subgroup from the class $\mathcal E$;
\item[(iii)] if a group from the class $\mathcal E$ acts on a tree,
it fixes a point, an end, or a pair of ends;
\item[(iv)] each maximal infinite subgroup from $\mathcal E$
is equal to its normalizer in $\pi$.
\end{description}
\end{dfn}
Note that the condition (iii) holds
if any group in $\mathcal E$ is amenable~\cite{Neb}.

\begin{thm}\label{accessability result}
Let $\pi$ be a torsion-free finitely presented group
and $\mathcal E$ be an elementary class of subgroups of $\pi$.
Then there exists an integer $K>0$ and a finite sequence 
$\pi_0,\pi_1,\dots ,\pi_m$ of subgroups of $\pi$ such that

$(1)$ $\pi_m=\pi$, and

$(2)$ for each $k$ with $0\le k<K$,
the group $\pi_k$ either belongs to $\mathcal E$ or
does not split as a nontrivial
amalgamated product or an HNN-extension over
a group from $\mathcal E$, and

$(3)$ for each $k$ with $K\le k\le m$, the group $\pi_k$ 
is the fundamental group of a finite graph of groups
with edge groups from $\mathcal E$, vertex groups from
$\{\pi_0,\pi_1,\dots ,\pi_{k-1}\}$,
and proper edge-to-vertex homomorphisms.
\end{thm}

\begin{prop}\label{virt nilp is elementary}
Let $\pi$ be a finitely presented group. 
If $\mathcal M_{a, b, \pi, n}\neq\emptyset$, then 
the class of virtually nilpotent subgroup of $\pi$
is elementary.
\end{prop}
\begin{proof} 
The proof is straightforward and can be found in~\cite{Bel5}.
\end{proof}
%
\begin{prop}\label{finiteness for vertices in graphs of groups} 
Let $\pi$ be the fundamental group of a finite 
graph of groups with virtually nilpotent edge groups.
Assume $\mathcal M_{a, b, \pi, n}\neq\emptyset$. Then

$(1)$ $\pi$ is finitely presented iff
all the vertex groups are finitely presented, and

$(2)$ $\dim_{\mathbb Q}H^{*}(\pi, \mathbb Q)<\infty$ iff
$\dim_{\mathbb Q} H^{*}(\pi_v, \mathbb Q)<\infty$
for every vertex group $\pi_v$.
%
\end{prop}
\begin{proof} Since $\mathcal M_{a, b, \pi, n}\neq\emptyset$,
any virtually nilpotent subgroup of $\pi$ is 
finitely generated~\cite{Bow2} and hence is 
the fundamental group
of a closed aspherical manifold~\cite{FH}.
The statement $(2)$ now follows from the Mayer-Vietoris 
sequence. For details on the proof of $(1)$ see~\cite{Bel5}.
\end{proof}

\section{Invariants of maps and actions}

This section is a condensed version of~\cite[sections 3;4]{Bel3}. 

\subsection{Invariants of continuous maps}

\begin{dfn} 
Let $B$ be a topological space and $S_B$ be a set.
Let $\iota$ be a map that, given a smooth manifold $N$, 
and a continuous map from $B$ into $N$,
produces an element of $S_B$. 
We call $\iota$ an {\it invariant of maps of $B$}
if the two following conditions hold:

$\mathrm{(1)}$ Homotopic maps $f_1:B\to N$ and 
$f_2\co B\to N$ have the same invariant.

$\mathrm{(2)}$ Let $h\co  N\to L$ be a
diffeomorphism of $N$ onto an open subset of $L$.
Then, for any continuous map $f\co B\to N$, 
the maps $f\co B\to N$ and $h\circ f\co B\to L$
have the same invariant.
\end{dfn}

There is a version of this definition for maps into
oriented manifolds.
Namely, we require that the target manifold is oriented and 
the diffeomorphism $h$ preserves orientation.
In that case we say that $\iota$ is an 
{\it invariant of maps into oriented manifolds}. 

\begin{ex}\bf \!(Tangent bundle.)\rm\
\label{Tangent bundle}
Assume $B$ is paracompact and $S_B$ is the set 
of isomorphism classes of real
vector bundles over $B$.
Given a continuous map 
$f\co B\to N$, set $\tau (f\co B\to N)=f^\#TN$,
the isomorphism class of
the pullback of the tangent bundle to $N$
under $f$.
Clearly, $\tau$ is an invariant.
\end{ex}

\begin{ex}\bf \!(Intersection number in oriented $n$-manifolds.)\rm\ 
\label{Intersection number in oriented $n$-manifolds}
Assume $B$ is compact
and fix two homology classes 
$\alpha\in H_{m}(B)$ and $\beta\in H_{n-m}(B)$.
(In this paper we always use singular (co)homology
with integer coefficients unless stated otherwise.)
Let $f\co B\to N$ be a continuous map of a
compact topological space $B$ into
an oriented $n$-manifold $N$ where $\mathrm{dim}(N)=n$.
Set $I_{n,\alpha,\beta}(f)$ to be the 
intersection number
of $f_*\alpha$ and $f_*\beta$ in $N$. 
It was verified in~\cite{Bel3} that $I_{n,\alpha,\beta}$
is an integer-valued invariant of maps into 
oriented manifolds.

\end{ex}
We say that an invariant of maps is {\it liftable}
if in the part $\mathrm{(2)}$ of the definition
the word ``diffeomorphism'' can be replaced by a
``covering map''.
For example, tangent bundle is a liftable invariant.
Intersection numbers are not liftable.
The following proposition shows to what extent 
it can be repaired.
 
\begin{prop}\rm\!~\cite{Bel3}\it\label{liftable prop} Let $p:\tilde{N}\to N$ be a
covering map of manifolds and let $B$ 
be a finite connected CW-complex.
Suppose that $f:B\to\tilde{N}$ is a 
map such that $p\circ f:B\to N$ is an embedding 
(i.e.,\! a homeomorphism onto its image).
Then $\iota(f)=\iota(p\circ f)$
for any invariant of maps $\iota$.
\end{prop}

\subsection{Invariants of actions}\label{Invariants of actions}

Assume $X$ is a smooth contractible manifold and 
let $\mathrm{Diffeo}(X)$ be the 
group of all self-diffeomorphisms of $X$ equipped with
the compact-open topology.
Let $\pi$ be the fundamental group of a 
finite-dimensional CW-complex $K$
with the universal cover $\tilde{K}$.

To any action $\rho\co \pi\to\text{Diffeo}(X)$,
we associate a 
continuous $\rho$-equivariant map $\tilde{K}\to X$ as follows.
Consider the $X$--bundle
$\tilde{K}\times_\rho X$ over $K$
where $\tilde{K}\times_\rho X$
is the quotient of $\tilde{K}\times X$
by the following action of $\pi$ 
$$\gamma(\tilde k,x)=(\gamma(\tilde k), \rho(\gamma)(x)),
\ \ \gamma\in\pi.$$ 
Since $X$ is contractible, the bundle has a section
that is unique up to homotopy through sections.
Any section can be lifted to
a $\rho$-equivariant continuous map 
$\tilde{K}\to \tilde{K}\times X$.
Projecting to $X$, we get a $\rho$-equivariant continuous map 
$\tilde{K}\to X$.
Note that any two $\rho$-equivariant continuous maps 
$\tilde{g}, \tilde{f}\co \tilde{K}\to X$, 
are $\rho$-equivariantly homotopic.
(Indeed, $\tilde{f}$ and $\tilde{g}$ descend to sections
$K\to \tilde{K}\times_\rho X$ that must be homotopic.
This homotopy lifts to a $\rho$--equivariant homotopy
of $\tilde{f}$ and $\tilde{g}$.) 
 
Assume now that $\rho(\pi)$ acts freely and 
properly discontinuously on $X$.
Then the map $\tilde{f}$ descends to a continuous map 
$f\co K\to X/\rho(\pi_1(K))$.
We say that $\rho$ is {\it induced} by $f$.

Let $\iota$ be an invariant of continuous maps of $K$.
Given an action $\rho$ such that
$\rho(\pi)$ acts freely and 
properly discontinuously on $X$, set
$\iota(\rho)$ to be $\iota(f)$
where $\rho$ is induced by $f$.
We say $\rho$ is an {\it invariant of free, proper discontinuous
actions of $\pi_1(K)$}. 
Similarly, any invariant $\iota$ of continuous maps of $K$
into oriented manifolds 
defines an invariant of free, proper discontinuous,
orientation-preserving actions on $X$.  

Note that actions conjugate 
by a diffeomorphism $\phi$ of $X$ have same invariants.
(Indeed, if $\tilde f\co\tilde{K}\to X$ is a 
$\rho$-equivariant map, the map $\phi\circ\tilde f$
is $\phi\circ\rho\circ\phi^{-1}$-equivariant.)
The same is true for invariants of orientation-preserving
actions when $\phi$ is orientation-preserving.

\begin{ex}\bf \!(Tangent bundle.)
\label{invariant tau}\rm\
Let $\tau$ be the invariant of maps defined in~\ref{Tangent bundle}.
Then, for any action $\rho$ such that
$\rho(\pi)$ acts freely and 
properly discontinuously on $X$, 
let $\tau(\rho)$ be the pullback of the tangent 
bundle to $X/\rho(\pi)$ via a map $f\co K\to X/\rho(\pi)$ 
that induces $\rho$.
\end{ex}

\begin{ex}\bf \!(Intersection number for
orientation preserving actions.)\rm\
Assume the cell complex $K$ is finite
and choose an orientation on $X$ (which
makes sense because, like any contractible
manifold, $X$ is orientable).
Given homology classes 
$\alpha\in H_{m}(K)$ and $\beta\in H_{n-m}(K)$,
let $I_{n,\alpha,\beta}$ is an invariant
of maps defined in~\ref{Intersection number in oriented $n$-manifolds}
where $n=\mathrm{dim}(X)$.

Let $\rho$ be an action of $\pi$
into the group of orientation preserving
diffeomorphisms of $X$ such that
$\rho(\pi)$ acts freely and 
properly discontinuously on $X$.
Then let $I_{n,\alpha,\beta}(\rho)$ be the 
intersection number
of $f_*\alpha$ and $f_*\beta$ in $X/\rho(\pi)$ 
where $f\co K\to X/\rho(\pi)$ is a map
that induces $\rho$.
\end{ex}

\subsection{Main theorem}

Throughout this section $K$ is a finite, 
connected CW-complex with a reference point $q$.
Let $\tilde{K}$ be the universal cover of $K$,
$\tilde{q}\in\tilde{K}$ be a preimage of $q\in K$.
Using the point $\tilde{q}$ we identify $\pi_1(K,q)$ with
the group of automorphisms of the covering $\tilde{K}\to K$.
Let $\iota$ is an invariant of actions.

\begin{thm}\label{main thm} 
Let $\rho_k:\pi=\pi_1(K, q)\to\mathrm{Isom}(X_k)$
be a sequence of isometric actions
of $\pi_1(K)$ on Hadamard $n$-manifolds $X_k$ such that
$\rho_k(\pi)$ is a discrete subgroup of $\mathrm{Isom}(X_k)$
that acts freely. 
Suppose that, for some $p_k\in X_k$,
$(X_k, p_k, \rho_k(\pi))$ converges in the
equivariant pointed
Lipschitz topology to $(X, p, \Gamma)$
and $(X_k, p_k, \rho_k)$ converges to $(X, p, \rho)$
in the pointwise convergence topology. 
Let $\tilde{h}: \tilde{K}\to X$ be a $\rho$-equivariant continuous map
and let $\bar{h}:K\to X/\Gamma$ be the drop of $\tilde h$.
Then $\iota(\rho_k)=\iota(\bar{h})$ for all large $k$.
\end{thm}

\begin{proof} Recall that according to~\ref{Invariants of actions}
such a map $\tilde h$ always exists and is unique up to
$\rho$-equivariant homotopy.

Let $\tilde{q}\in F\subset\tilde{K}$ be a finite subcomplex
that projects onto $K$.
Clearly, the finite set 
$S=\{\gamma\in\pi_1(K,q): \gamma(F)\cap F\neq\emptyset\}$
generates $\pi_1(K,q)$.

Note that $X$ is contractible. (Indeed,
any spheroid in $X$ lies in the diffeomorphic
image of a metric ball in $X_j$.
Any metric ball in a Hadamard manifold is contractible.
Thus $\pi_*(X)=1$.) 

Choose $\epsilon >0$ so small that 
$\tilde{h}(F)$ lies in the open
ball $B(p,1/10\epsilon)\subset X$.
For large $k$, we find an $\epsilon$-Lipschitz
approximation $(\tilde{f}_k, \tilde{g}_k, \phi_k, \tau_k)$ 
between $(X_k,p_k,\rho_k(\pi))$ and
$(X,p,\Gamma)$.  

By lemma~\ref{approximation} we can assume that
$\tau_k(\rho(\gamma ))=\rho_k(\gamma)$ for all $\gamma\in S$.
Hence, the map
$\tilde{h}_k=\tilde{g}_k\circ\tilde{h}:F\to X_k$
is $\rho_k$-equivariant.
Extend it by equivariance to a $\rho_k$-equivariant
map $\tilde{h}_k:\tilde{K}\to X_k$.
Passing to quotients we get a map
$h_k:K\to X_k/\rho_k(\pi_1(K))$
such that $\iota(\rho_k)=\iota(h_k)$.

By construction, $h_k=g_k\circ \bar{h}$
where $g_k$ is the drop of $\tilde{g}_k$.
Since $g_k$ is a diffeomorphism,
$\iota(h_k)=\iota(\bar{h})$.
Hence $\iota(\rho_k)=\iota(\bar{h})$ and
the proof is complete. 
\end{proof}

\begin{rmk} There is a version of the theorem~\ref{main thm} 
for invariants of maps into oriented manifolds.
Suppose all the actions $\rho_k$
on $X_k$ preserve orientations  
(it makes sense because
being a contractible manifold $X_k$ is orientable).
Fix an orientation on $X$ (which is also
contractible) and choose orientations of $X_k$
so that diffeomorphisms $g_k$ preserve orientations.
Then the same proof gives $\iota(\rho_k)=\iota(\bar{h})$
for any invariant of maps into oriented manifolds $\iota$.
Yet this new orientation on $X_k$ may be different from 
the original one.
\end{rmk}

\begin{cor}\label{when iota(rho.n)=iota(rho)}
Suppose that in addition to the assumptions of the 
theorem\ \rm\ref{main thm}\it\
any of the following holds
\begin{itemize}
\item $\iota$ is a liftable invariant, or
\item $\rho(\pi_1(K))=\Gamma$, or
\item $\bar{h}$ is homotopic to an embedding.
\end{itemize}
Then $\iota(\rho_k)=\iota(\rho)$ for all large $k$.
\end{cor}
\begin{proof} 
Let $h:K\to X/\rho(\pi)$ be the drop of $\tilde h$.
Since $\iota(\rho)=\iota(h)$,
it suffices to understand when $\iota(h)=\iota(\bar{h})$.
This is trivially true if $\iota$ is liftable or
if $\rho(\pi_1(K))=\Gamma$.
In case $\bar{h}$ is homotopic to an embedding,~\ref{liftable prop}
implies that $\iota(h)=\iota(\bar{h})$.
\end{proof}
\begin{rmk}
In Kleinian group theory the condition 
$\rho(\pi_1(K))=\Gamma$ means, by definition,
that $\rho_k$ converges to $\rho$ {\it strongly}. 
\end{rmk}

\section{Applications}

\subsection{Strongly pinched manifolds have zero Pontrjagin classes}

Let $X$ be the limit of Hadamard $n$-manifolds
$X_k$ in the pointed Lipschitz topology.
Assume that the sectional curvature of $X_k$
is within $[-1-1/k, -1]$.
Then $X$ is a smooth $n$-manifold~\cite{Fuk}
and by a standard argument $X$ is an inner metric space
of curvature $\ge -1$ and $\le -1$
in the sense of Alexandrov.
Hence~\cite{Ale} implies that
$X$ is isometric to the real hyperbolic space.

\begin{thm}
Let $\pi$ be the fundamental group of a finite
aspherical complex.
Suppose that $\pi$ is not virtually nilpotent
and that $\pi$ does not split as a 
nontrivial amalgamated product
or an HNN-extension over a virtually nilpotent group.
Then, for any positive integer $n$, there exists
an $\epsilon=\epsilon(\pi, n)>0$ such that any manifold in
the class $\mathcal M_{-1-\epsilon, -1, \pi, n}$
is tangentially homotopy equivalent to
a manifold of constant negative curvature.
\end{thm}
\begin{proof}
Arguing by contradiction, consider a sequence
of manifolds $X_k/\rho_k(\pi)\in\mathcal M_{-1-1/k, -1,\pi,n}$
that are not tangentially homotopy equivalent to
manifolds of constant negative curvature. 
  
According to~\ref{divergence implies split},
$(X_k, p_k,\rho_k)$ is precompact in the pointwise convergence topology,
for some $p_k\in X_k$.
Hence~\ref{pointwise implies lipsch} implies that
$(X_k, p_k,\rho_k(\pi))$ is precompact in the 
equivariant pointed Lipschitz topology.
Pass to subsequence so that $(X_k, p_k,\rho_k)$
converges to $(X, p,\rho)$ in the pointwise convergence topology
and $(X_k, p_k,\rho_k(\pi))$
converges to $(X, p,\Gamma)$ in the 
equivariant pointed Lipschitz topology.

By~\ref{geom contains alg, etc}, 
$\rho$ is a free action of a discrete group $\pi$, in particular
$X/\rho(\pi)$ is a manifold with fundamental 
group isomorphic to $\pi$.
Set $K=X/\rho(\pi)$ and let $\tau$ be the invariant of 
representations of $\pi_1(K)$ defined in~\ref{invariant tau}. 
According to~\ref{when iota(rho.n)=iota(rho)}, 
$\tau(\rho)=\tau(\rho_k)$ for large $k$,
thus $\rho_k\circ\rho^{-1}$
induces a tangential homotopy
equivalence of $X_k/\rho_k(\pi)$ and $X/\rho(\pi)$. 
Since $X$ is isometric to the real hyperbolic space, 
the proof is complete.
\end{proof}

\begin{thm}\label{no split implies zero pontrj}
Let $\pi$ be the fundamental group of a
finite-dimensional aspherical CW-complex $K$ such that
$\dim_\mathbb Q\oplus_m H^{4m}(K,\mathbb Q)<\infty$. 
Suppose that $\pi$ is finitely presented,
not virtually nilpotent
and that $\pi$ does not split as a 
nontrivial amalgamated product
or an HNN-extension over a virtually nilpotent group.
Then, for any positive integer $n$, there exists
an $\epsilon=\epsilon(\pi, n)>0$ such that 
the tangent bundle of any manifold in
the class $\mathcal M_{-1-\epsilon, -1, \pi, n}$
has zero rational Pontrjagin classes.
\end{thm}
\begin{proof}
Arguing by contradiction, consider a sequence
of manifolds $X_k/\rho_k(\pi)\in\mathcal M_{-1-1/k, -1,\pi,n}$
whose tangent bundles have nonzero rational Pontrjagin classes.
According to~\ref{divergence implies split},
$(X_k, p_k,\rho_k)$ is precompact in the pointwise convergence topology,
for some $p_k\in X_k$.
Hence~\ref{pointwise implies lipsch} implies that
$(X_k, p_k,\rho_k(\pi))$ is precompact in the 
equivariant pointed Lipschitz topology.
Pass to subsequence so that $(X_k, p_k,\rho_k)$
converges to $(X, p,\rho)$ in the pointwise convergence topology
and $(X_k, p_k,\rho_k(\pi))$
converges to $(X, p,\Gamma)$ in the 
equivariant pointed Lipschitz topology.

By an elementary homological argument
there exists a finite connected subcomplex
$X\subset K$ such that the inclusion $i:X\to K$
induces a $\pi_1$--epimorphism
and a $\oplus_m H^{4m}(-,\mathbb Q)$-monomorphism
(for example, see~\cite{Bel5} for a proof).
The sequence of isometric actions 
$\rho_k\circ i_*$ of $\pi_1(X,x)$ on $X_k$ satisfies
the assumptions of~\ref{main thm}, therefore,
$\tau(\rho_k\circ i_*)=\tau(\rho\circ i_*)$ for all large $k$.
In other words, the pullback bundles $i^\#\tau(\rho_k)$
are isomorphic to the vector bundle $i^\#\tau(\rho)$.

By~\ref{geom contains alg, etc}, 
$\rho$ is a free action of a discrete group $\pi$, in particular
$X/\rho(\pi)$ is a real hyperbolic manifold.
Hence, the tangent bundle to $X/\rho(\pi)$ has zero
rational Pontrjagin classes $p_m$ for $m>0$. 
(Avez observed in~\cite{Ave} that
Pontrjagin forms on any conformally flat 
manifold vanish.)
Hence $\tau(\rho)$ and $i^\#\tau(\rho)$ have zero
rational Pontrjagin classes because they are pullbacks
of the tangent bundle to $X/\rho(\pi)$.
Thus $i^\#\tau(\rho_k)$ has zero rational Pontrjagin classes
for all large $k$. 
Hence $i^*p_m(\tau(\rho_k))=p_m(i^\#\tau(\rho_k))=0$
all large $k$.
Since $i^*$ is injective, $p_m(\tau(\rho_k))=0$ 
which is a contradiction. 
\end{proof}

\begin{rmk} For any connected finite-dimensional cell complex $B$,
the Pontrjagin character defines an isomorphism
$\widetilde{KO}(B)\to\oplus_{m>0}H^{4m}(B,\mathbb Q)$.
In particular, a vector bundle $\xi$
has zero rational Pontrjagin classes iff
for some $n$ the Whitney sum 
$\underbrace{\xi\oplus\xi\oplus\cdots\oplus\xi}_{n}$
is a trivial bundle.
\end{rmk}

\begin{prop}\label{cohom mayer-vietoris}
Let $\pi$ be the fundamental group of a finite graph of groups
such that each edge group has cohomological dimension $\le 2$.
Let $n$ be a positive integer. 
Assume that there exists
an $\epsilon>0$ such that, for any
vertex group $\pi_v$,
the tangent bundle of any manifold in
the class $\mathcal M_{-1-\epsilon, -1, \pi_v, n}$
has zero rational Pontrjagin classes.
Then the tangent bundle of any manifold in
the class $\mathcal M_{-1-\epsilon, -1, \pi, n}$
has zero rational Pontrjagin classes.
\end{prop}
\begin{proof}
Following~\cite{SW}, we assemble the cell complexes
$K(\pi_v,1)$ and $K(\pi_e,1)\times [-1,1]$
into an $K(\pi,1)$ cell complex by using edge--to--vertex
monomorphisms.
By Mayer-Vietoris sequence the map
$H^{4m}(\pi,\mathbb Q)\to\oplus_v H^{4m}(\pi_v,\mathbb Q)$
induced by inclusions $\pi_v\to\pi$ is
an isomorphism.
Let $N\in\mathcal M_{-1-\epsilon, -1, \pi, n}$ 
and let $N_v$ be the cover of $N$ that corresponds
the inclusions $\pi_v\to\pi$; clearly
$N_v\in\mathcal M_{-1-\epsilon, -1, \pi_v, n}$. 
Since coverings preserve tangent bundles,
the map $H^{4m}(N,\mathbb Q)\to H^{4m}(N_v,\mathbb Q)$
takes $p_m(TN)$ to $p_m(TN_v)$.
Thus, $p_m(TN)=0$ iff $p_m(TN_v)=0$ for all $v$. 
\end{proof}

\begin{thm}\label{small nil subg imply zero pontrj} 
Let $\pi$ be a finitely presented group 
with finite $4k$th Betti numbers for all $k$
Assume that any nilpotent
subgroup of $\pi$ has cohomological dimension $\le 2$.
Then for any $n$ there exist $\epsilon=\epsilon(n,\pi)>0$ 
such that the tangent bundle of any manifold in the class
$\mathcal M_{-1-\epsilon, -1, \pi, n}$ 
has zero rational Pontrjagin classes.
\end{thm}
\begin{proof}
Applying the theorem~\ref{accessability result}, 
we get a sequence $\pi_0,\pi_1,\dots ,\pi_m$ 
of subgroups of $\pi$.
In particular, for every $k$ with $0\le k<K$, the group $\pi_k$ 
either is virtually nilpotent or does not split over
a virtually nilpotent subgroup of $\pi$.
By the proposition~\ref{finiteness for vertices in graphs of groups}, 
the group $\pi_k$ is finitely presented and  
has finite Betti numbers.
Therefore, by~\ref{no split implies zero pontrj},
any manifold in the class 
$\mathcal M_{-1-\epsilon, -1, \pi_k, n}$ 
has zero rational Pontrjagin classes.

Repeatedly applying~\ref{cohom mayer-vietoris}, 
we deduce that any manifold in the class 
$\mathcal M_{-1-\epsilon, -1, \pi, n}$ 
has zero rational Pontrjagin classes.
\end{proof}

\begin{cor}\label{pinch hyper gr}
Let $\pi$ be a word-hyperbolic group.
Then for any $n$ there exist $\epsilon=\epsilon(n,\pi)>0$ 
such that the tangent bundle of any manifold in the class
$\mathcal M_{-1-\epsilon, -1, \pi, n}$ 
has zero rational Pontrjagin classes.
\end{cor}
\begin{proof}
Any torsion free word-hyperbolic group is the fundamental group of
a finite aspherical cell complex~\cite[5.24]{CDP}.
Moreover, any virtually nilpotent subgroup of a 
torsion free word-hyperbolic group
of is either trivial or infinite cyclic.
So~\ref{small nil subg imply zero pontrj} applies.
\end{proof}

\subsection{Intersections in negatively
curved manifolds}

Given an oriented $n$-manifold $N$,
consider the intersection form $I_{m,n-m}(\cdot,\cdot)$
of type $(m, n-m)$
$$I_{m,n-m}:H_m(N)\otimes H_{n-m}(N)\to \mathbb Z.$$

A homotopy equivalence $f$ of orientable $n$-manifolds is called
{\it $(m, n-m)$-intersection preserving} if, for some choice of orientations
of the manifolds,
$$I_{m,n-m}(f_*\alpha,f_*\beta)=I_{m,n-m}(\alpha,\beta)$$
for all $(\alpha,\beta )\in H_m(N)\otimes H_{n-m}(N)$.

We now extend the notion of $(m, n-m)$-intersection preserving
homotopy equivalence to non-orientable manifolds.
Recall that a homotopy equivalence of manifolds $f:N\to L$
is called {\it orientation-true} if the vector bundles
$TN$ and $f^\# TL$ have equal first Stiefel-Whitney class.
Let $f$ be an orientation-true homotopy equivalence of
non-orientable manifolds. 
Then $f$ lifts to a homotopy equivalence
$\tilde{f}:\tilde{N}\to \tilde{L}$ of orientable
two-fold covers.
We say that $f$ is {\it $(m, n-m)$-intersection preserving} if 
so is $\tilde{f}$.

If $f$ is {\it $(m, n-m)$-intersection preserving} for all $m$,
we say that $f$ is {\it intersection preserving}.
For example, if $f$ is homotopic to a 
homeomorphism, then $f$ is $(m, n-m)$-intersection preserving
for all $m$ [Dold, 13.21].

\begin{thm}\label{bound on intersections}
Let $K$ be a connected finite cell complex and
let $I_{n,\alpha,\beta}$ be the invariant of maps defined
in\ \rm\ref{Intersection number in oriented $n$-manifolds}.\it\
Let $\rho_k:\pi_1(K)\to\mathrm{Isom}(X_k)$
be a sequence of free, isometric actions
of $\pi_1(K)$ on Hadamard $n$-manifolds $X_k$. 
Suppose that, for some $p_k\in X_k$,
$(X_k, p_k, \rho_k(\pi_1(K)))$ is precompact in the
equivariant pointed Lipschitz topology and 
$(X_k, p_k, \rho_k)$  is precompact in 
the pointwise convergence topology.

Then, for any sequence of continuous maps 
$f_k:K\to X_k/\rho_k(\pi_1(K))$ that induce $\rho_k$,
the sequence of integers $I_{n,\alpha,\beta}(f_k)$
is bounded.
\end{thm}

\begin{proof} 
Since $I_{n,\alpha,\beta}(f_k)=I_{n,\alpha,\beta}(\rho_k)$,  
this is a particular case of~\ref{main thm}.
\end{proof}

\begin{thm}\label{int pres: main result}
Let $\pi$ finitely generated group
and let $\rho_k:\pi\to\mathrm{Isom}(X_k)$
be a sequence of free, isometric actions
of $\pi$ on Hadamard $n$-manifolds $X_k$. 
Suppose that, for some $p_k\in X_k$,
$(X_k, p_k, \rho_k(\pi))$ is precompact in the
equivariant pointed Lipschitz topology and 
$(X_k, p_k, \rho_k)$  is precompact in 
the pointwise convergence topology. 
Assume that the Betti numbers 
$b_m$ and $b_{n-m}$ of $\pi$ are finite.
Then the set of manifolds $\{X_k/\rho_k(\pi)\}$
falls into finitely many $(m, n-m)$-intersection preserving
homotopy types.
\end{thm}

\begin{proof} Argue by contradiction.
We can assume that no two manifold from  $\{X_k/\rho_k(\pi)\}$
are $(m, n-m)$-intersection preserving
homotopy equivalent. 
Set $K=X_1/\rho_1(\pi)$.
Let $f_k:K\to X_k/\rho_k(\pi)$
be a homotopy equivalence that induces $\rho_k$.
We next show that, after passing to a subsequence,
the homotopy equivalence $f_{k+1}\circ f_k^{-1}$
is $(m, n-m)$-intersection preserving which yields
a contradiction.

Since $\pi_1(K)$ is finitely generated, 
the group $H^1(K,\mathbb Z_2)\cong\mathrm{Hom}(\pi_1(K),\mathbb Z_2)$
is finite. So passing to subsequence we
can assume that the vector bundles $f_k^\#TX_k/\rho_k(\pi)$
have the same first Stiefel-Whitney class $w$.
Hence the homotopy equivalences $f_{k+1}\circ f_k^{-1}$
are orientation-true.

The first Stiefel-Whitney class $w$ defines a two-fold-cover
$\tilde{K}\to K$ and an index two subgroup $\tilde\pi=\pi_1(\tilde{K})$
of $\pi_1(K)$. Restricting $\rho_k$ to $\tilde{\pi}$ 
we get a sequence of free, isometric actions
of $\tilde{\pi}$ on Hadamard $n$-manifolds $X_k$. 
Notice that $(X_k, p_k, \rho_k(\tilde{\pi}))$ is precompact in the
equivariant pointed Lipschitz topology and 
$(X_k, p_k, \rho_k|_{\tilde{\pi}})$  is precompact in 
the pointwise convergence topology.
Let $\tilde{f}_k:\tilde{K}\to X_k/\rho_k(\tilde{\pi})$
be a homotopy equivalence that induces $\rho_k|_{\tilde{\pi}}$.

Take arbitrary 
$\alpha\in H_m(\tilde{K})$ and $\beta\in H_{n-m}(\tilde{K})$
and let $I_k$ be the intersection number
of $\tilde{f}_{k*}\alpha$ and $\tilde{f}_{k*}\beta$ in $X_k/\rho_k(\tilde{\pi})$.
By an elementary homological argument
there exists a finite connected subcomplex $K^\prime$
of $\tilde{K}$ such that the inclusion $i:K^\prime\subset K$
induces epimorphisms of the fundamental groups
and $i_*H_*(K^\prime)$ contains $\alpha$ and $\beta$
(for example, see~\cite{Bel5} for a proof).

Let $i_*\alpha^\prime=\alpha$ and $i_*\beta^\prime=\beta$.
Clearly, $I_k$ is equal to the intersection number
of $\tilde{f}_{k*}i_*\alpha^\prime$ and $\tilde{f}_{k*}i_*\beta^\prime$ in 
$X_k/\rho_k(\tilde{\pi})$.
According to~\ref{bound on intersections}, the set of integers $\{I_k\}$ is finite.
Hence, passing to subsequence, we can assume that the intersection number
of $\tilde{f}_{k*}\alpha$ and $\tilde{f}_{k*}\beta$ in $X_k/\rho_k(\tilde{\pi})$
is independent of $k$.

The groups $H_m(\tilde{K})$ and $H_{n-m}(\tilde{K})$ have finite rank
since the Betti numbers $b_m$ and $b_{n-m}$ of $\tilde{K}$
are finite. 
(Recall that an abelian group $A$ has finite rank
if there exists a finite subset $S\subset A$
such that any non-torsion element of 
$A$ is a linear combination of elements of $S$.)
Therefore, the intersection form is determined by
intersection numbers of finitely many 
homology classes. 
(Torsion elements do not matter because
the intersection number of a torsion class
and any other class is zero.)
Then the argument of the previous paragraph implies that,
passing to subsequence,
we can assume that the intersection number
of $\tilde{f}_{k*}\alpha$ and $\tilde{f}_{k*}\beta$ in $X_k/\rho_k(\tilde{\pi})$
is independent of $k$ for any classes $\alpha$ and $\beta$.
In other words, $\tilde{f}_{k+1}\circ\tilde{f}_k^{-1}$
is an $(m, n-m)$-intersection preserving homotopy equivalence.
\end{proof}

\begin{cor}
Let $\pi$ be a finitely generated group with finite 
Betti numbers.
Assume $\rho_k\co\pi\to\mathrm{Isom}(X_k)$
is a sequence of free, isometric actions
of $\pi$ on Hadamard $n$-manifolds $X_k$. 
Suppose that, for some $p_k\in X_k$,
$(X_k, p_k, \rho_k(\pi))$ is precompact in the
equivariant pointed Lipschitz topology and 
$(X_k, p_k, \rho_k)$  is precompact in 
the pointwise convergence topology. 
Then, for each $m$, the set of manifolds $\{X_k/\rho_k(\pi)\}$
falls into finitely many $(m, n-m)$-intersection preserving
homotopy types.\qed
\end{cor}

\begin{prop}\label{int pres:nilpotent groups}
Let $N$ be an orientable pinched negatively curved manifold with
virtually nilpotent fundamental group. Then the
intersection number of any two homology
classes in $N$ is zero.
\end{prop}
\begin{proof}
It suffices to prove that $N$ is homeomorphic
to $\mathbb R\times Y$ for some space $Y$.
First note that any torsion free, discrete, 
virtually nilpotent group $\Gamma$
acting on a Hadamard manifolds $X$ of pinched
negative curvature must have either one or two 
fixed points at infinity~\cite[3.3.1]{Bow2}. 
If $\Gamma$ has only one fixed point,
$\Gamma$ is parabolic and, hence, it preserves
all horospheres at the fixed point.
Therefore, if $H$ is such a horosphere,
$X/\Gamma$ is homeomorphic to $\mathbb R\times H/\Gamma$.
If $\Gamma$ has two 
fixed points, $\Gamma$ preserves a bi-infinite geodesic.
Hence $X/\Gamma$ is the total space of a vector bundle over
a circle. Thus $X/\Gamma$ is homeomorphic
to $\mathbb R\times Y$ for some space $Y$ unless $X/\Gamma$
is the M\"obius band which is impossible since
$N=X/\Gamma$ is orientable.
\end{proof}

\begin{cor}\label{int pres: finiteness for groups that do not split}
Assume that $\pi$ is a finitely presented 
group does not split over a virtually nilpotent group.
Let $m\le n$ be integers such that the Betti numbers 
$b_m$ and $b_{n-m}$ of $\pi$ are finite.
Then, for any $a\le b<0$,
the class $\mathcal M_{a, b, \pi, n}$ breaks into finitely
many $(m, n-m)$-intersection preserving
homotopy types.
\end{cor}
\begin{proof} By~\ref{int pres:nilpotent groups} we
can assume that $\pi$ is not virtually nilpotent.
Let $N_k$ be an arbitrary sequence of manifolds
from $\mathcal M_{a, b, \pi, n}$ 
represented as $X_k/\rho_k(\pi)$.
According to~\ref{divergence implies split},
$(X_k, p_k,\rho_k)$ is precompact in the pointwise convergence topology,
for some $p_k\in X_k$.
Hence~\ref{pointwise implies lipsch} implies that
$(X_k, p_k,\rho_k(\pi))$ is precompact in the 
equivariant pointed Lipschitz topology and 
we are done because of~\ref{int pres: main result}.
\end{proof}

\begin{thm}\label{no split means bound on intersections}
Let $K$ be a connected finite cell complex such
that the group $\pi_1(K)$ does not split 
over a virtually nilpotent group. 
Let $I_{n,\alpha,\beta}$ be the invariant of maps defined
in\ \rm\ref{Intersection number in oriented $n$-manifolds}\it\ 
and let $a\le b<0$ be real numbers.
Then, for any sequence of continuous maps $f_k:K\to N_k$ that 
induce isomorphisms of fundamental groups of $K$
and $N_k\in\mathcal M_{a,b,\pi_1(K),n}$,
the sequence of integers $I_{n,\alpha,\beta}(f_k)$
is bounded.
\end{thm}

\begin{proof} By~\ref{int pres:nilpotent groups} we
can assume that $\pi$ is not virtually nilpotent.
Arguing by contradiction, consider a sequence
of maps $f_k\co K\to N_k$ that induce $\pi_1$-isomorphisms
of $K$ and manifolds $N_k\in\mathcal M_{a,b,\pi_1(K),n}$
and such that no two integers $I_{n,\alpha,\beta}(f_k)$
are equal. 

Each map $f_k$ induces a free, properly 
discontinuous action $\rho_k$ of $\pi_1(K)$ 
on the universal cover $X_k$ of $N_k$
such that $I_{n,\alpha,\beta}(f_k)=I_{n,\alpha,\beta}(\rho_k)$.  
According to~\ref{divergence implies split},
$(X_k, p_k,\rho_k)$ is precompact in the pointwise convergence topology,
for some $p_k\in X_k$.
Hence~\ref{pointwise implies lipsch} implies that
$(X_k, p_k,\rho_k(\pi))$ is precompact in the 
equivariant pointed Lipschitz topology and 
we are done because of~\ref{bound on intersections}.
\end{proof}

\subsection{Vector bundles with negatively curved total spaces}

The following is a slight generalization of the 
theorem~\ref{intro:vector bundles}. The proof makes
use of the classification theorem proved
in the appendix.

\begin{thm}\label{embeddings finite}
Let $M$ be a closed negatively curved
manifold of dimension $\ge 3$.
Suppose $n>\dim(M)$ is an integer and $a\le b<0$ are real numbers. 
Let $f_k:M\to N_k$ be a sequence of smooth embeddings of $M$ into
manifolds $N_k$ such that for each $k$
\begin{itemize}
\item $f_k$ induces a monomorphism of fundamental groups, and
\item $N_k$ is a complete Riemannian $n$-manifold 
with sectional\\ curvatures within $[a,b]$. 
\end{itemize}
Then the set of 
the normal bundles $\nu (f_k)$ of the embeddings 
falls into finitely many isomorphism classes.
In particular, up to diffeomorphism, 
only finitely many manifolds from the class
$\mathcal M_{a, b, \pi_1(M), n}$ are total spaces of
vector bundles over $M$.
\end{thm}
\begin{proof} Passing to covers corresponding to $f_{k*}$,
we can assume that $f_k$ induce isomorphisms of fundamental groups.
Arguing by contradiction, assume that $\nu_k=\nu (f_k)$ are
pairwise nonisomorphic. 

First, reduce to the case when the bundles $\nu_k$ are orientable. 
Since $H^1(M,\mathbb Z_2)$ is finite,
we can pass to subsequence and assume that
the bundles $\nu_k$ have equal first Stiefel-Whitney classes.
Pass to two-fold covers corresponding to this
Stiefel-Whitney class. 
Then $f_k$ lift to embeddings with orientable
normal bundles. Using~\ref{classifying nonorientable bundles},
we can pass to subsequence so that these normal bundles
are pairwise nonisomorphic.
Furthermore, a finite cover of $M$ is still a closed negatively curved
manifold of dimension $\ge 3$ and any finite cover of $N_k$
is a complete Riemannian manifold with
sectional curvatures within $[a,b]$.
Thus, it suffices to consider the case of orientable $\nu_k$.

Each map $f_k$ induces a free, properly 
discontinuous action $\rho_k$ of $\pi_1(K)$ 
on the universal cover $X_k$ of $N_k$
such that $\tau(f_k)=\tau(\rho_k)$.  
According to~\ref{divergence implies split},
$(X_k, p_k,\rho_k)$ is precompact in the pointwise convergence topology,
for some $p_k\in X_k$.
Hence~\ref{pointwise implies lipsch} implies that
$(X_k, p_k,\rho_k(\pi))$ is precompact in the 
equivariant pointed Lipschitz topology. 
According to~\ref{main thm}, 
the set of vector bundles $\{\tau(f_k)\}$
falls into finitely many isomorphism classes.

In particular, there are only finitely many possibilities for the
total Pontrjagin class of $\tau(f_k)$.
The normal bundle $\nu_k$ of the embedding $f_k$ satisfies
$\nu_k\oplus TM\cong\tau(f_k)$. Applying the total
Pontrjagin class, we get $p(\nu_k)\cup p(TM)=p(\tau(f_k))$.
The total Pontrjagin class of any bundle is a unit,
hence we can solve for $p(\nu_k)$.
Thus, there are only finitely many possibilities for $p(\nu_k)$.
 
Thus, according to~\ref{classifying orientable bundles},
we can pass to subsequence so that the (rational) Euler classes 
of $\nu_k$ are all different.
Denote the integral Euler class by $e(\nu_k)$.

First, assume that $M$ is orientable. 
Recall that, by definition, the Euler class
$e(\nu_k)$ is the image of the 
Thom class $\tau(\nu_k)\in H^m(N_k,N_k\backslash f_k(M))$
under the map $f_{k}^*\co H^m(N_k,N_k\backslash f_k(M))\to H^m(M)$.
According
to~\cite[VIII.11.18]{Dol} the Thom class   
has the property 
$\tau(\nu_k)\cap [N_k,N_k\backslash f_k(M)]=f_{k*}[M]$
where $[N_k,N_k\backslash f_k(M)]$ is 
the fundamental class of the pair 
$(N_k,N_k\backslash f_k(M))$ and $[M]$ is
the fundamental class of $M$.
Therefore, for any $\alpha\in H_m(M)$,
the intersection number of
$f_{k*}\alpha$ and $f_{k*}[M]$ in $N_k$ satisfies
$$I(f_{k*}[M], f_{k*}\alpha)=
\langle\tau(\nu_k), f_{k*}\alpha\rangle =
\langle f^*\tau(\nu_k), \alpha\rangle =
\langle e(\nu_k),\alpha\rangle .$$

Since $M$ is compact, $H_m(M)$ is finitely generated;
we fix a finite set of generators.
The (rational) Euler classes are all different, hence
the homomorphisms 
$\langle e(\nu (f_k)),-\rangle\in\text{Hom}(H_m(M),\mathbb Z)$ 
are all different. 
Then there exists a generator
$\alpha\in H_m(M)$ such that 
$\{\langle e(\nu (f_k)),\alpha\rangle\}$ is an infinite set of integers.
Hence $\{I(f_{k*}[M], f_{k*}\alpha)\}$ is an infinite set
of integers. Combining~\ref{no split means bound on intersections} 
and~\ref{aspherical do not split}, 
we get a contradiction.

Assume now that $M$ is nonorientable.
Let $q\co\tilde M\to M$ be the orientable
two-fold cover. 
Any finite cover of aspherical manifolds
induces an injection on rational
cohomology~\cite[III.9.5(b)]{Bro}.
Hence $e(q^\#\nu (f_k))=q^*e(\nu (f_k))$
implies that the rational Euler classes of 
the pullback bundles $q^\#\nu (f_k)$ are all different, and
there are only finitely many possibilities for
the total Pontrjagin classes of $q^\#\nu (f_k)$.

Furthermore, the bundle map $q^\#\nu (f_k)\to\nu (f_k)$
induces a smooth two-fold cover of the total spaces,
thus the total space of $q^\#\nu (f_k)$ belongs to
$\mathcal M_{a, b, \pi_1(\tilde{M}), n}$.
and we get a contradiction 
as in the oriented case.
\end{proof}

\begin{rmk}
More generally, the theorem~\ref{embeddings finite} is true whenever
$M$ is a closed smooth aspherical manifold such that
no finite index subgroup of $\pi_1(M)$ splits
over a virtually nilpotent group.
The proof we gave works {\it verbatim}. 
In particular, we can take $M$ to be a smooth manifold that
satisfies the conditions of~\ref{aspherical do not split}.
\end{rmk}

\begin{rmk}\label{vector bundles infinite}
In some cases it is easy to decide when
there exist {\it infinitely}
many vector bundles of the same rank over a given base.
Namely, it suffices to check that
certain Betti numbers of the base are nonzero.
For example, by a simple $K$-theoretic argument
the set of isomorphism classes of
rank $m$ vector bundles over a finite cell complex 
$B$ is infinite provided $m\ge\dim (B)$
and $\oplus_m H^{4m}(B,\mathbb Q)\neq 0$.
In fact, any element of $\oplus_m H^{4m}(B,\mathbb Q)$
is the Pontrjagin character of some vector bundle over $B$.
Furthermore, 
oriented rank two vector bundles over $B$ are in one-to-one
correspondence with $H^2(B,\mathbb Z)$
via the Euler class. 

Note that many arithmetic closed real hyperbolic manifolds
have nonzero Betti numbers in all dimensions~\cite{MR}.
Any closed complex hyperbolic manifold has nonzero
even Betti numbers because the powers
of the K\"ahler form are noncohomologous to zero.
Similarly, for each $k$, closed quaternion hyperbolic
manifolds have nonzero $4k$th Betti numbers.
\end{rmk}

\section{Digression: pinching invariants}

We now reinterpret our results 
in terms of certain natural pinching invariants.
See insightful discussions in~\cite{Gro2, Gro3} 
for more information.

For a smooth manifold $N$ that admits a metric of 
pinched negative curvature, define $pinch(N)\in [1,\infty )$
to be the infimum of the numbers $a/b$ such that
$N$ is diffeomorphic to a manifold in
the class $\mathcal M_{a, b, \pi_1(N), \dim(N)}$.

\begin{ex}
Replacing ``diffeomorphic'' by ``homeomorphic''
leads to a different invariant.
For all $n\ge 6$,  Farrell, Jones, and Ontaneda constructed sequences of 
closed negatively curved $n$-manifolds $N_k$ with $pinch(N_k)\to 1$ 
such that each $N_k$ is homeomorphic but not diffeomorphic
to a hyperbolic manifold~\cite{FJO}.
\end{ex}

\begin{ex} For any $n\ge 4$,
Gromov and Thurston~\cite{GT} constructed 
a sequence of closed $n$-manifolds $N_k$
such that $pinch(N_k)>1$ and $pinch(N_k)$ converges to $1$
as $k\to\infty$. These manifolds are not diffeomorphic
to closed real hyperbolic manifolds.
Furthermore, they constructed a
sequence of closed negatively curved $n$-manifolds $N_k$
such that $pinch(N_k)\to\infty$.
\end{ex}
\begin{ex}
Let $\xi_k$ be a sequence of pairwise non-isomorphic rank $m$
vector bundles over a closed
negatively curved manifold $M$ of dimension $\ge 3$. 
Then according to the theorem~\ref{intro:vector bundles}, 
their total spaces $E(\xi_k)$ have the property that 
$pinch(N_k)$ converges to infinity.
\end{ex}
\begin{ex}
Suppose that $N$ is a manifold of pinched negative
curvature with word-hyperbolic fundamental group
and nontrivial total Pontrjagin class. 
Then theorem~\ref{pinch hyper gr} implies $pinch(N)>1$.
\end{ex}

Let $\pi$ be the fundamental group
of a manifold of pinched negative curvature. 
Define $pinch_n(\pi)\in [1,\infty )$ to be 
the infimum over the numbers $a/b$
such that $\mathcal M_{a, b, \pi, n}\neq\emptyset$.
Clearly, $pinch(N)\ge pinch_{\dim(N)}(\pi_1(N))$.
Notice that $$pinch_n(\pi)\ge pinch_{n+1}(\pi)$$
because if $N\in\mathcal M_{a, b, \pi, n}$, then there is a
warped product metric on $N\times\mathbb R$ such that
$N\times\mathbb R\in\mathcal M_{a, b, \pi, n+1}$~\cite{FJ5}. 
Clearly, if $\Gamma$ is a subgroup of $\pi$, then
$pinch_n(\Gamma)\le pinch_n(\pi)$.

\begin{ex}
Let $\pi$ be a finitely generated group
that does not act on an $\mathbb R$-tree with no global
fixed point and virtually nilpotent arc stabilizers. 
Then according to the propositions~\ref{pointwise implies lipsch}
and~\ref{divergence implies split}, there always
exists a smooth $n$-manifold $N$ with complete $C^{1,\alpha}$
Riemannian metric of bounded Alexandrov curvature
with pinching equal to $pinch_n(\pi)$. 
Moreover, by a result of Nikolaev~\cite{Nik}, this metric on $N$
can be approximated in (nonpointed) Lipschitz topology
by complete Riemannian metrics on $N$
whose pinchings converge
to $pinch_n(\pi)$. In particular,
$pinch(N)=pinch_n(\pi)$.

Furthermore, if $pinch_n(\pi)=1$, 
then by~\ref{no split plus pinched means hyperbolic} the manifold
$N$ carries a complete real hyperbolic metric. 
Thus, for a group $\pi$ as above, 
either $\pi$ is the fundamental group of
a real hyperbolic manifold, or else $pinch_n(\pi)>1$.
\end{ex}
\begin{ex}
Any closed negatively curved $n$-manifold $N$
with $pinch_n(\pi_1(N))=1$ 
is diffeomorphic to a real hyperbolic manifold.
This is obvious if $\dim(N)=2$ 
and follows from~\ref{no split plus pinched means hyperbolic}
and~\ref{aspherical do not split} if $\dim(N)>2$.
\end{ex}
%
%
\begin{ex} Let $\pi$ be a (discrete)
group with Kazhdan's property $(T)$.
Then $\pi$ is finitely generated and
any action of $\pi$ on an $\mathbb R$-tree has a global
fixed point~\cite{dlHV}. Furthermore, $\pi$ is not the fundamental
group of a real hyperbolic manifold~\cite{dlHV}.
Thus we conclude that $pinch_n(\pi)>1$.
\end{ex}
Sometimes it is possible to compute or at least estimate
$pinch_n(\pi)$. Here we only give two examples
that use harmonic maps. See~\cite{Gro3} for other
results in this direction.
\begin{ex}
Let $M$ be a closed K\"ahler manifold  
such that $\pi_1(M)$ has cohomological dimension $>2$.
Then $pinch_n(\pi_1(M))\ge 4$ for all $n\ge\dim(M)$~\cite{YZ}.
If, in addition, $M$ is complex hyperbolic, then
$pinch_n(\pi_1(M))=4$.
\end{ex} 
\begin{ex} Let $M$ be a closed quaternion hyperbolic or
Cayley hyperbolic manifold.
Then $pinch_n(\pi_1(M))=4$. Moreover,
if $\pi$ is a quotient of $\pi_1(M)$, 
then $pinch_n(\pi)\ge 4$~\cite{MSY}.
\end{ex} 

\appendix
\section{Classifying vector bundles}

The purpose of this appendix is to prove that
the isomorphism type of any vector bundle 
is determined, up to finitely many possibilities,
by the characteristic classes of the bundle.
This fact is apparently well-known to experts,
yet there seem to be no published proof. 
The proof given below is elementary and mainly
uses obstruction theory.
I am most grateful to Jonathan Rosenberg from whom I
learned the orientable case and to Sergei P.~Novikov
for help in the non-orientable case. 
Here is the precise statement for orientable vector bundles;
the proof can be found in~\cite{Bel3}. 

\begin{thm}\label{classifying orientable bundles}
Let $K$ be a finite CW-complex and $m$ be a positive integer.
Then the set of isomorphism classes of oriented real 
(complex, respectively) rank $m$ vector bundles over $K$
with the same rational Pontrjagin classes and the rational Euler class 
(rational Chern classes, respectively) is finite. 
\end{thm}
First, we review the proof for orientable 
vector bundles of, say, even rank $m$. 
Characteristic classes can be thought of as homotopy classes 
of maps from the classifying space $BSO(m)$ to Eilenberg-MacLane spaces. 
For example, Euler class and Pontrjagin classes are given by
$e\in H^{m}(BSO(m), \mathbb Z)\cong [BSO(m), K(m, \mathbb Z)]$
and $p_i\in H^{4i}(BSO(m), \mathbb Z)\cong [BSO(m), K(4i, \mathbb Z)]$.
The map $(e, p_1,\dots, p_{m/2-1})$ of $BSO(m)$
to the product of Eilenberg-MacLane spaces is known to induce an 
isomorphism in rational cohomology. 
Thus, since the spaces are simply-connected,
$(e, p_1,\dots, p_{m/2-1})$ is a rational homotopy equivalence. 
Now the obstruction theory implies that a map of a 
finite cell complex $K$ into the product of Eilenberg-MacLane spaces
can have only finitely many nonhomotopic liftings to $BSO(m)$.
In other words, only finitely many vector bundles over $K$
can have the same characteristic classes.

The above argument  
fails for nonorientable bundles, 
due to the fact $BO(m)$ is not simply connected 
(i.e.,\! the map $c=(p_1,\dots, p_{[m/2]})$
of $BO(m)$ to the product of Eilenberg-MacLane spaces
is {\it not} a rational homotopy equivalence
even though it induces an isomorphism on
rational cohomology). 
Yet essentially the same result is true. 
To make a precise statement we need the
following background.

Given a finite CW-complex $K$,
the set of isomorphism classes of
rank $m$ vector bundles over a $K$ is in one-to-one correspondence
with the set of homotopy classes of maps $[K,BO(m)]$;
to a map $f:K\to BO(m)$ there corresponds the pullback 
$f^\#\gamma_m$ of the universal rank $m$ vector bundle 
$\gamma_m$ over $BO(m)$.
A vector bundle $f^\#\gamma_m$ is orientable iff $f$ lifts to 
$BSO(m)$ which is a two-fold-cover of $BO(m)$.
In terms of characteristic classes a vector bundle  
is orientable iff its first Stiefel-Whitney class
vanishes.

Let $\xi=f^\#\gamma_m$ and $\eta=g^\#\gamma_m$ be nonorientable
vector bundles that have the same first Stiefel-Whitney class 
$w_1(\xi)=w_1(\eta)=w\in H^1(K,\mathbb Z_2)\setminus 0$.
The universal coefficient theorem  provides
a natural isomorphism of $H^1(K,\mathbb Z_2)$
and 
$\mathrm{Hom}(H_1(K),\mathbb Z_2)\cong\mathrm{Hom}(\pi_1(K),\mathbb Z_2)$.
Thus, to a nonzero element of $H^1(K,\mathbb Z_2)$ there corresponds an
epimorphism of $\pi_1(K)$ onto $\mathbb Z_2$ whose kernel is an index two
subgroup of $\pi_1(K)$. This index two
subgroup defines a two-fold-cover $\tilde{K}\to K$.

Let $p:\tilde{K}\to K$ be the two-fold-cover that corresponds
to the class $w$. 
Clearly $p^*w=0$, hence the pullback bundles $p^\#\xi$
and  $p^\#\eta$ are orientable.
In other words, the maps $f\circ p$ and $g\circ p$
of $\tilde{K}$ to $BO(m)$ can be lifted to $BSO(m)$. 
The lifts $\tilde{f},\tilde{g}:\tilde{K}\to BSO(m)$
are equivariant with respect the covering actions of $\mathbb Z_2$.
Clearly, $f$ is homotopic to $g$ iff $\tilde{f}$ is equivariantly
homotopic to $\tilde{g}$. 

\begin{thm}\label{classifying nonorientable bundles}
Let $\xi$ be nonorientable vector bundle over a 
finite CW-complex $K$ and let $p:\tilde{K}\to K$ be 
the two-fold-cover that corresponds
to the first Stiefel-Whitney class of $\xi$.
Then $\xi$ is 
is determined up to finitely many possibilities
by the Euler class and the total Pontrjagin class of $p^\#\xi$.
\end{thm}

The proof is based on the equivariant obstruction theory
which is reviewed below.
%
%
Suppose $\tilde{K}\to K$ is a two-fold-cover
of a finite CW-complex $K$. 
Thus, we get an involution on the set of cells of $\tilde{K}$
and hence an involution $\iota$ of the cellular chain complex
$C_*(\tilde{K})$. Clearly, $\iota$ commutes with the boundary
homomorphism.  

The ``usual'' cellular cohomology $H^*(K,\Pi)$
of $\tilde{K}$
with coefficients in a finitely generated abelian group $\Pi$
is the homology of the complex 
$\mathrm{Hom}(C_*(\tilde{K}),\Pi)$.

Let $\mathbb Z_2$ act (by group automorphisms) 
on $\Pi$; in other words there is $i\in\mathrm{Aut}(\Pi)$
such that $i^2$ is the identity.
A cochain $f\in\mathrm{Hom}(C_*(\tilde{K}),\Pi)$ is called
{\it equivariant} if $f(\iota c)=if(c)$.
The {\it equivariant} cellular cohomology 
$H^*_e(\tilde{K},\Pi)$ of $\tilde{K}$
with coefficients in $\Pi$ is the homology of the subcomplex
of equivariant cochains
$$\mathrm{Hom}_{\mathbb Z_2}(C_*(\tilde{K}),\Pi)
\subset\mathrm{Hom}(C_*(\tilde{K}),\Pi).$$ 
This inclusion induces a homomorphism 
$H^*_e(\tilde{K},\Pi)\to H^*(\tilde{K},\Pi)$. 

\begin{lem}\label{finite kernel} The homomorphism 
$H^*_e(\tilde{K},\Pi)\to H^*(\tilde{K},\Pi)$
has finite kernel.
\end{lem}
\begin{proof} Denote the complex 
$\mathrm{Hom}(C_*(\tilde{K}),\Pi)$ by $C$
and the subcomplex of equivariant cochains 
$\mathrm{Hom}_{\mathbb Z_2}(C_*(\tilde{K}),\Pi)$
by $C_+$.
Let $C_-\subset C$ be the subcomplex
of anti-equivariant cochains where a cochain $f$
is called
{\it anti-equivariant} if $$f(\iota c)=-if(c).$$

First notice that the inclusion $C_+\to C_+\oplus C_-$
induces a monomorphism in homology.
Indeed, take $f\in C_+$ such that $(f,0)$ is a boundary,
that is $(f,0)=\partial(g, h)=(\partial g, \partial h)$.
Thus, $f=\partial g$.

Second, show that the map $C_+\oplus C_-\to C$
that takes $(f,g)$ to $f+g$ has finite kernel.
Indeed, assume $f+g=0$. 
Since $f$ is equivariant and $g$ is anti-equivariant
we deduce
$$if(c)=f(\iota c)=-g(\iota c)=ig(c).$$
Hence $f=g$ and, therefore, both $f$ and $g$ have
order two. 
Thus, the kernel of $C_+\oplus C_-\to C$ lies in the $2$-torsion
of $C_+\oplus C_-$.
In particular, the kernel lies in the
torsion subgroup of $C\oplus C$ which is finite
because $C$ is a finitely generated abelian group.
(In fact, if $\tilde{K}$ has $k$ cells,
$C_*(\tilde{K})$ is a free abelian group of rank $k$.
Hence $C$ is a the direct sum of $k$ copies of a finitely
generated group $\Pi$.)
 
Third, prove that the map $C_+\oplus C_-\to C$
has finite cokernel. Notice that its image $C_++ C_-$ 
contains $2C$. (Indeed $2f$ is the sum of
an equivariant cochain $f_+$ and an anti-equivariant cochain
$f_-$ where $f_\pm(c)$ is defined as $f(c)\pm f(\iota c)$.) 
Thus, it suffices to check that $C/2C$ is finite
which is true because $C$ is a finitely generated abelian group.

Finally, two short exact sequences of chain complexes
$$0\to\mathrm{ker}\to C_+\oplus C_-\to C_++ C\to 0\ \ \
\mathrm{and}\ \ \
0\to C_++ C_-\to C\to \mathrm{coker}\to 0$$
induce long exact sequences in homology with $H(\mathrm{ker})$
and $H(\mathrm{coker})$ finite.
Therefore, both $H(C_+\oplus C_-)\to H(C_++ C_-)$ and
$H(C_++ C_-)\to H(C)$ have finite kernel.
In particular, the composition 
$$H(C_+)\to H(C_+\oplus C_-)\to H(C_++ C_-)\to H(C)$$
has finite kernel as desired.
\end{proof}
%
%
For free, proper discontinuous actions
the usual obstruction theory routinely generalizes
to the equivariant case~\cite{Dug}.

Let $\tilde{K}\to K$ and $\tilde{Y}\to Y$ be two-fold-covers
of CW-complexes and let $\tilde{f}$ and $\tilde{g}$ be continuous maps
of $\tilde{K}$ into $\tilde{Y}$ that are equivariant with respect
to a unique isomorphism of covering groups.
Assume that $\tilde{Y}$ is simply-connected.
Then $\tilde{f}$ and $\tilde{g}$ are equivariantly homotopic on
the one-skeleton~\cite[Lemma~9.1]{Dug}.
Assume $\tilde{f}$ and $\tilde{g}$ are equivariantly homotopic on
the $(n-1)$-skeleton.

Consider the difference cochain $d^n(\tilde{f},\tilde{g})$ 
with coefficients in $\pi_n(\tilde Y)$ that
comes from the usual obstruction theory; in fact, $d^n(\tilde{f},\tilde{g})$
is a cocycle since $\tilde{f}$ and $\tilde{g}$ are defined
on the whole $\tilde{K}$.
It was shown in~\cite[p266]{Dug} that $d^n(\tilde{f},\tilde{g})$
is an equivariant cochain and, furthermore, 
if $d^n(\tilde{f},\tilde{g})$ is an equivariant coboundary,
then $\tilde{f}$ and $\tilde{g}$
are equivariantly homotopic on
the $n$-skeleton. 
   
In our case $Y=BO(m)$ and $\tilde Y=BSO(m)$. 
In particular, $\pi_n(BO(m))$ is a finitely generated
abelian group hence the lemma~\ref{finite kernel} applies.

\begin{proof}[Proof of the theorem~\ref{classifying nonorientable bundles}]
We are going to argue by contradiction.
Let $\xi_i$ be an infinite sequence of vector bundles
given by maps $f_i:K\to BO(m)$. Assume that
the bundles $\xi_i$ have equal first Stiefel-Whitney
classes and let $p:\tilde{K}\to K$ be the two-fold-cover
corresponding to this first Stiefel-Whitney
class. Assume also that the bundles
$p^\#\xi_i$ have the same Euler class
and total Pontrjagin class.

As before, let $\tilde{f}_i:\tilde{K}\to BSO(m)$
be the lift of $f_i$.
Note that all the maps $\tilde{f}_i$  
are equivariantly homotopic on
the one-skeleton~\cite[Lemma~9.1]{Dug}.
Using the theorem~\ref{classifying orientable bundles}, 
we pass to subsequence so that
all the bundles $p^\#\xi_i$ are isomorphic.
In other words, all the maps $\tilde{f}_i$ are
(non-equivariantly) homotopic.

Let $n>1$ be the smallest integer such that
infinitely many of $\tilde{f}_i$'s are not equivariantly
homotopic on the $n$-skeleton. 
Pass to subsequence to assume that
$\tilde{f}_i$'s are equivariantly
homotopic on the $(n-1)$-skeleton.
Thus, the difference cochains $d^n(\tilde{f}_i,\tilde{f}_1)$
are defined.

Since all the maps $\tilde{f}_i$ are
(non-equivariantly) homotopic, $d^n(\tilde{f}_i,\tilde{f}_1)$
represents the zero element in the (non-equivariant)
cohomology.
By lemma~\ref{finite kernel}, we can pass to subsequence so that
the difference cochains $d^n(\tilde{f}_i,\tilde{f}_1)$
represents the same element in the equivariant
cohomology. 

By the properties of the difference cochain
$d^n(\tilde{f}_i,\tilde{f}_j)=d^n(\tilde{f}_i,\tilde{f}_1)-
d^n(\tilde{f}_j,\tilde{f}_1)$, hence
$d^n(\tilde{f}_i,\tilde{f}_j)$
represents the zero element in the
equivariant cohomology.
Hence, $f_i$ and $f_j$ are equivariantly homotopic on the
$n$-skeleton. This is a contradiction
with the assumption that the sequence $\{f_i\}$
is infinite.\end{proof}

\begin{rmk} \label{pontrjagin classifies}
Note that, if either
$m$ is odd or $m>\mathrm{dim}(K)$, the rational Euler class is
zero, and, hence rational Pontrjagin classes determine
a vector bundle up to a finite number of possibilities.
(We use here that $H^1(K,\mathbb Z_2)$ is finite.) 
\end{rmk}

\small
\bibliographystyle{amsalpha}
\bibliography{vbxxx}
\end{document}